\documentclass[11pt]{amsart}

\usepackage{amsmath}
\usepackage{amssymb}
\usepackage{mathtools}
\usepackage{mathrsfs}
\usepackage{dsfont}

\usepackage{graphicx}
\usepackage{tikz}
\usetikzlibrary{arrows.meta,positioning}

\usepackage{microtype}

\usepackage[colorlinks=true,
            linkcolor=blue,
            citecolor=blue,
            urlcolor=blue]{hyperref}

\usepackage[nameinlink,noabbrev]{cleveref}

\newtheorem{theorem}{Theorem}[section]
\newtheorem{proposition}[theorem]{Proposition}
\newtheorem{lemma}[theorem]{Lemma}
\newtheorem{corollary}[theorem]{Corollary}

\theoremstyle{definition}

\theoremstyle{remark}

\numberwithin{equation}{section}
\usepackage[margin=1in]{geometry}
\title{Torsion Points in Families of Elliptic Curves}

\author{Ester Sgallová}

\address{Department of Mathematics,
Michigan State University,
East Lansing, MI 48824, USA}

\email{sgallova@msu.edu}

\date{}



\begin{document}

\begin{abstract}
We study the specialization of torsion points in one-parameter degenerations of elliptic curves with singular fibers of Kodaira type. We describe the resulting torsion curves, their ramification,
intersection multiplicities, and irreducible and connected components.
\end{abstract}

\maketitle
\section{Introduction}
Let $\pi\colon\mathcal E\to\Delta$ be a relatively minimal elliptic
fibration over a complex disk, with zero section $Z$. We assume that
$\mathcal E$ is smooth, that $\pi$ is smooth over
$\Delta^*:=\Delta\setminus\{0\}$, and that the central fiber
$\mathcal E_0:=\pi^{-1}(0)$ is a singular Kodaira fiber. We use the
Kodaira classification and the standard descriptions of local monodromy
\cite{Kodaira1963,Miranda1989,SilvermanAdvanced}.

For $t\in\Delta^*$, write $E_t=\pi^{-1}(t)$
and let $E_t[n]$ denote its group of $n$-torsion
points, for $n\geq1$. As $t$ varies over
$\Delta^*$, these groups form a finite étale
cover of degree $n^2$. We define the
$n$-torsion curve $T_n$ as the Zariski
closure of this cover in $\mathcal E$:
\[
    T_n:=
    \overline{\bigcup_{t\in\Delta^*}E_t[n]}
    \subset\mathcal E.
\]

Fix $t\in\Delta^*$. Monodromy around $0$
acts on $E_t[n]\simeq(\mathbb Z/n\mathbb Z)^2$,
and its orbits correspond to the irreducible
components of $T_n$. An orbit
$O\subset E_t[n]$ of cardinality $r$
determines an irreducible component
$B_O\subset T_n$ whose normalization
$\widetilde B_O\to\Delta$ has degree $r$
and is totally ramified above $0$.
In suitable analytic coordinates, this map
is given by $t=u^r$. In particular,
$B_O$ is a section of $\pi$ if and only
if $r=1$.

Let $p_O\in\mathcal E_0$ denote the unique
specialization point of $B_O$.
Monodromy determines the irreducible
components of $T_n$ and their degrees
over $\Delta$, but not their specialization
points. These may lie in the smooth locus
of the reduced central fiber or at its
singular points. Moreover, distinct
irreducible components may specialize
to the same point and thus belong to
the same connected component of $T_n$. The purpose of this paper is to determine this specialization geometry.

We distinguish three classes of singular Kodaira fibers.
Fibers of type $I_k$, with $k\geq 1$, have
\emph{multiplicative reduction}.
Fibers of type $I_k^*$, with $k\geq 1$, have
\emph{additive, potentially multiplicative reduction}:
after a finite ramified base change, they acquire
multiplicative reduction.
The remaining singular Kodaira fibers
\[
    II,\ III,\ IV,\ I_0^*,\ IV^*,\ III^*,\ II^*
\]
have \emph{additive, potentially good reduction}:
after a finite ramified base change, they acquire
good reduction.

For each of these fiber types, we determine the
specialization points of the irreducible components
of $T_n$, their intersection multiplicities with
the irreducible components of the central fiber,
and the numbers of irreducible and connected
components of $T_n$.

There are several settings in which such torsion curves occur. On an
elliptically fibered K3 surface, irreducible components of torsion
multisections provide examples of constant cycle curves
\cite[Section~6.2]{Huybrechts2014}. They also appear in the work of Bogomolov and
Tschinkel on potential density of rational points on elliptic K3
surfaces, where degrees and genera of multisections are used to
distinguish torsion from non-torsion multisections \cite[Section~3.4]{BogomolovTschinkel2000}.
A second motivation comes from the boundary of toroidal
compactifications of Shimura varieties. In ongoing joint work on
$\mathcal A_2$, we compare generating series of Humbert surfaces and
Noether--Lefschetz loci via a Shimura lift and study the closures of
these cycles in a toroidal compactification. The resulting boundary
calculations involve torsion curves on modular elliptic surfaces.

Another motivation comes from the study of isogenies
in families. Recent work of Lin establishes a finiteness result for height values in geometric isogeny classes of motives under suitable hypotheses, including an adelic Mumford–Tate condition
\cite{LinIsogenyHeights}.
For elliptic curves, a cyclic isogeny of degree
$n$ is determined by a cyclic subgroup of order
$n$ in $E_t[n]$. Over $\Delta^*$, the cover
parametrizing such subgroups is obtained from
the locus of torsion points of exact order $n$
by identifying points that generate the same
subgroup. The specialization and intersection
calculations for $T_n$ provide local geometric
information that may be useful for studying
degenerations of these covers and related
questions about isogenies in families.

\medskip
\noindent\textbf{Main results.}
We begin with fibers of type $I_k$. Write
\[
    \mathcal E_0=C_0+\cdots+C_{k-1},
\]
where the components are labeled cyclically and
the zero section meets $C_0$. The indices on
the $C_j$ are taken modulo $k$.

Let $P_1,\ldots,P_k$ be the nodes of the
central fiber in cyclic order, where $P_j$
lies between the branches labeled $C_{j-1}$
and $C_j$. For $k=1$, both branches at $P_1$
belong to $C_0$.

For $0\leq j\leq k-1$, set
\[
    C_j^*
    :=C_j\setminus\operatorname{Sing}(\mathcal E_0).
\]
For a zero-dimensional intersection, we use
$\operatorname{length}$ for its total length and
$\operatorname{length}_p$ for its local length at a point $p$.

\begin{theorem}[Specialization over an $I_k$-fiber]\label{thm:intro-Ik}
Let $k\geq1$ and set
\[
  g:=\gcd(n,k),\qquad m:=\frac{n}{g},\qquad
  d_i:=\gcd(ik,n)\quad(0\leq i\leq n-1).
\]
Then the following statements hold.
\begin{enumerate}
\item For $0\leq j\leq k-1$,
\[
  \operatorname{length}(T_n\cap C_j^*)
  =\begin{cases}
      n,& k\mid jn,\\
      0,& k\nmid jn.
    \end{cases}
\]
\item For $1\leq j\leq k$, the local intersection length is
\[
  \operatorname{length}_{P_j}(T_n\cap\mathcal E_0)
  =n\,\#\{1\leq i\leq n-1:(j-1)n<ik<jn\}.
\]
For each $i$ in this set there are $d_i$ irreducible components
$B_{i,c}$, indexed by $c\in\mathbb Z/d_i\mathbb Z$, specializing to
$P_j$. Each has one analytic branch at $P_j$, and the normalization
of $B_{i,c}$ has degree $n/d_i$ over $\Delta$.
\item If $k\geq2$ and $(j-1)n<ik<jn$, then
\[
  B_{i,c}\cdot C_{j-1}=\frac{jn-ik}{d_i},\qquad
  B_{i,c}\cdot C_j=\frac{ik-(j-1)n}{d_i},
\]
and $B_{i,c}$ has intersection number zero with the other
irreducible components of the central fiber. If $ik=jn$, then
$d_i=n$; the corresponding $B_{i,c}$ are sections specializing to
$C_j^*$, each with $B_{i,c}\cdot C_j=1$.
\item The number of irreducible components of $T_n$ is
\[
  \#\operatorname{Irr}(T_n)
  =g^2\sum_{d\mid m}d\,\varphi\!\left(\frac{m}{d}\right),
\]
where $\varphi$ is Euler's totient function.
\item The number of connected components of $T_n$ is
\[
  \#\pi_0(T_n)=gn+\min\{n-g,k\}.
\]
\item For $0\leq j\leq k-1$, the intersection number with $C_j$ is
\[
  T_n\cdot C_j
  =\sum_{r\in\mathbb Z}\bigl(n-|nj-rk|\bigr)_+,
  \qquad (x)_+:=\max\{x,0\}.
\]
\item Write $n=\alpha k+\beta$, with
$0\leq\beta<k$. Then, for every
$0\leq j\leq k-1$,
\[
    T_n\cdot C_j
    =
    T_{\beta}\cdot C_j
    +\frac{n^2-\beta^2}{k}.
\]
Here we set $T_{0}\cdot C_j=0$.
\end{enumerate}
\end{theorem}

The difference between parts~(4) and~(5) arises because distinct irreducible
components may specialize to the same node and become connected there.
Parts~(6) and~(7) give the intersection numbers
of $T_n$ with the irreducible components of
the central fiber. In particular, part~(7)
shows that replacing $n$ by $n+k$ increases
all these intersection numbers by the same amount.

For fibers of type $I_k^*$ with $k\geq1$, label the multiplicity-two
components of the central fiber by $D_0,\ldots,D_k$ and their
successive intersection points by
\[
  P_j=D_{j-1}\cap D_j,\qquad 1\leq j\leq k.
\]
There are four further components, each of multiplicity one, two
meeting $D_0$ and two meeting $D_k$. With a suitable basis for
$H_1(\mathcal E_t,\mathbb Z)$, the local monodromy is $-M_k$, where
\[
  M_k=\begin{pmatrix}1&k\\0&1\end{pmatrix}.
\]
After the quadratic base change $t=s^2$, normalization and minimal
resolution give a semistable model with fiber of type $I_{2k}$.
The involution induced by $s\mapsto-s$ acts on the components of that
fiber by $C_j\mapsto C_{2k-j}$, after a suitable cyclic labeling.
The $I_k^*$-configuration is recovered by taking the quotient,
resolving its singularities, and passing to the relatively minimal
model. We use the specialization calculation on the $I_{2k}$-model
and this action of the involution to identify the specialization
points on the $I_k^*$-fiber.

For a statement \(A\), let
\(\mathbf 1_A\) denote its indicator.

\begin{theorem}[Specialization over an $I_k^*$-fiber]
\label{thm:intro-Ikstar}
Let $k\geq1$ and put $g:=\gcd(n,2k)$. Then:
\begin{enumerate}
\item No irreducible component of $T_n$ specializes to any of the
four nodes meeting the multiplicity-one components. For
$1\leq j\leq k$, the local intersection length in $P_j$ is
\[
  \operatorname{length}_{P_j}(T_n\cap\mathcal E_0)
  =2n\,\#\{1\leq i\leq n-1:(j-1)n<2ki<jn\}.
\]
For each $i$ in this set, there are $d_i:=\gcd(2ki,n)$
irreducible components $B_{i,c}$, indexed
by $c\in\mathbb Z/d_i\mathbb Z$, specializing
to $P_j$. The normalization of each $B_{i,c}$
is totally ramified over $0$, with ramification
index $2n/d_i$. For each such component $B_{i,c}$,
\[
  B_{i,c}\cdot D_{j-1}=\frac{jn-2ki}{d_i},\qquad
  B_{i,c}\cdot D_j=\frac{2ki-(j-1)n}{d_i}.
\]
\item The number of irreducible components is
\[
  \#\operatorname{Irr}(T_n)
   =\frac12\sum_{d\mid n}\varphi(d)\frac nd\gcd(d,2k)
      +\frac{\gcd(n,2)}2
      +\mathbf 1_{2\mid n}\mathbf 1_{4\mid nk}.
\]
\item The number of connected components is
\[
  \#\pi_0(T_n)
   =\frac{ng+\gcd(n,2)}2
    +\mathbf 1_{2\mid n}\mathbf 1_{4\mid nk}
    +\min\!\left\{k,\frac{n-g}{2}\right\}.
\]
\end{enumerate}
\end{theorem}

The specialization points away from the nodes, and the
intersection numbers with each irreducible component of the central
fiber, are given in Section~4. In particular, the formulas distinguish
the multiplicity-one outer components from the multiplicity-two chain.

Finally, consider fibers with potentially good reduction. Their
local monodromy has order $6$, $4$, $3$, or $2$ for types
$II$, $III$, $IV$, and $I_0^*$, respectively; the starred duals
$II^*$, $III^*$, and $IV^*$ have inverse monodromy, up to conjugacy.
After a ramified base change of degree equal to the monodromy order,
followed by normalization, resolution and relative minimalization,
the resulting family has smooth central fiber. The covering group
acts on this smooth elliptic curve and on its $n$-torsion. Comparing
its orbits with the quotient model determines the specialization
points on the original Kodaira fiber.

\begin{theorem}[Potentially good fibers]\label{thm:intro-good}
The fiber types in each pair
\[
  (II,II^*),\qquad (III,III^*),\qquad (IV,IV^*)
\]
have the same monodromy orbits on $ E_t[n]$, up to the
identification of their monodromy actions with a matrix and its
inverse. Consequently, for each pair they have the same number and
degrees of irreducible components of $T_n$.

For a fiber of type $I_0^*$,
\[
  \#\operatorname{Irr}(T_n)
   =\begin{cases}
      \dfrac{n^2-1}{2}+1, & 2\nmid n,\\[3pt]
      \dfrac{n^2-4}{2}+4, & 2\mid n.
    \end{cases}
\]
For a fiber of type $IV$ or $IV^*$,
\[
  \#\operatorname{Irr}(T_n)
   =\begin{cases}
      \dfrac{n^2-3}{3}+3, & 3\mid n,\\[3pt]
      \dfrac{n^2-1}{3}+1, & 3\nmid n.
    \end{cases}
\]
For a fiber of type $III$ or $III^*$,
\[
  \#\operatorname{Irr}(T_n)
   =\begin{cases}
      \dfrac{n^2-4}{4}+3, & 2\mid n,\\[3pt]
      \dfrac{n^2-1}{4}+1, & 2\nmid n.
    \end{cases}
\]
For a fiber of type $II$ or $II^*$,
\[
  \#\operatorname{Irr}(T_n)
   =\begin{cases}
      \dfrac{n^2-1}{6}+1,
       & 2\nmid n\text{ and }3\nmid n,\\[3pt]
      \dfrac{n^2-3}{6}+2,
       & 3\mid n\text{ and }2\nmid n,\\[3pt]
      \dfrac{n^2-4}{6}+2,
       & 2\mid n\text{ and }3\nmid n,\\[3pt]
      \dfrac{n^2-6}{6}+3,
       & 6\mid n.
    \end{cases}
\]

For types $II$, $III$, and $IV$, the irreducible components
corresponding to non-fixed monodromy orbits all specialize to the
unique singular point of the central fiber. They therefore belong
to one connected component of $T_n$, whereas the section components
are disjoint from them and from one another. For types
$I_0^*$, $IV^*$, $III^*$, and $II^*$, every specialization point lies
in the smooth locus of the reduced central fiber, and distinct
monodromy orbits give distinct specialization points. An orbit of
cardinality $r$ specializes to a component of the central fiber of
multiplicity $r$. Hence
\[
  \#\pi_0(T_n)=\#\operatorname{Irr}(T_n)
\]
for these four types.
\end{theorem}

The proofs are organized by reduction type. Section~2 identifies
monodromy orbits with the irreducible components of $T_n$.
Section~3 treats $I_k$-fibers using ramified base change and the
intersection theory of the resolved model. Section~4 treats
$I_k^*$-fibers via quadratic semistable reduction and the involution
of the $I_{2k}$-model. Sections~5 and~6 construct good-reduction
models for the remaining fiber types and determine the specialization
of their irreducible components.

\section{Monodromy and Components of the Torsion Curve}
\label{sec:monodromy-components}

We first describe the irreducible components of $T_n$ in terms of
monodromy over the punctured disk. Put
\[
    T_n^\circ:=T_n\cap\pi^{-1}(\Delta^*).
\]
By the definition of $T_n$, the fiber of
$T_n^\circ$ over $t\in\Delta^*$ is $E_t[n]$.
Since multiplication by $n$ on a smooth
elliptic family is finite étale in
characteristic zero, the map
\[
    T_n^\circ\longrightarrow\Delta^*
\]
is a finite unramified holomorphic covering
of degree $n^2$.

Fix $t\in\Delta^*$. Since
$T_n^\circ\to\Delta^*$ is a covering map,
the unique path-lifting property
\cite[\S~4]{ForsterRiemannSurfaces}
defines a monodromy action
\[
    \rho_n:\pi_1(\Delta^*,t)
    \longrightarrow\operatorname{Aut}(E_t[n]).
\]
The punctured disk has fundamental group
$\pi_1(\Delta^*,t)\simeq\mathbb Z$,
generated by a positively oriented loop
around the origin.

To identify this action, write $E_t\simeq\mathbb C/\Lambda_t$,
where $\Lambda_t$ is a lattice in $\mathbb C$. Then
\[
    E_t[n]
    \simeq \frac{\frac1n\Lambda_t}{\Lambda_t}
    \simeq \frac{\Lambda_t}{n\Lambda_t}
    \simeq \frac{H_1(E_t,\mathbb Z)}{nH_1(E_t,\mathbb Z)}.
\]
These identifications are compatible with transport in the family.
Consequently, the action of a loop on $E_t[n]$ is obtained by
reducing its action on $H_1(E_t,\mathbb Z)$ modulo $n$.
In particular, the local monodromy matrix on $n$-torsion is
the reduction modulo $n$ of the integral local monodromy matrix.

\begin{proposition}[Monodromy orbits and irreducible components]
\label{prop:monodromy-components}
The monodromy orbits on $E_t[n]$ are in bijection with the
irreducible components of $T_n$.

Let $O\subset E_t[n]$ be an orbit of cardinality $r$, and let
$B_O\subset T_n$ be the corresponding irreducible component.
Then the normalization
\[
    \widetilde B_O\longrightarrow\Delta
\]
has degree $r$ and a unique point above $0$. At this point the
map is given, in suitable analytic coordinates, by
\(
    t=u^r.
\)
In particular, its ramification index is $r$, and $B_O$ is a
section of $\pi$ if and only if $r=1$.
\end{proposition}

\begin{proof}
The connected components of a covering space correspond to the
orbits of the fundamental group on a fiber; see
\cite[\S\S~3--5]{ForsterRiemannSurfaces}.
Applied to $T_n^\circ\to\Delta^*$, this gives a bijection
between the monodromy orbits on $E_t[n]$ and the connected
components of $T_n^\circ$.

Each connected component of $T_n^\circ$ is a smooth connected
complex curve and hence is irreducible. Its closure in
$\mathcal E$ is an irreducible component of $T_n$. Conversely,
every irreducible component of $T_n$ meets $T_n^\circ$,
since $T_n$ is the Zariski closure of the relative torsion
locus over $\Delta^*$. Thus the monodromy orbits correspond
to the irreducible components of $T_n$.

Let $O$ be a monodromy orbit of cardinality $r$.
The corresponding connected component of $T_n^\circ$
is a degree-$r$ holomorphic covering of $\Delta^*$.
By the classification of finite holomorphic coverings
of the punctured disk
\cite[Thm.~5.10]{ForsterRiemannSurfaces},
this covering is isomorphic to
\[
    \Delta_u^*\longrightarrow\Delta_t^*,
    \qquad t=u^r.
\]

The restriction $B_O\to\Delta$ is proper, since
$\pi$ is proper. As $B_O$ is an irreducible curve
dominating $\Delta$, its fibers are zero-dimensional
and hence finite. Thus $B_O\to\Delta$ is a finite map.

The normalization $\widetilde B_O$ is a Riemann
surface, and its restriction over $\Delta^*$ is
the covering described above. Its normalized
extension over $0$ is therefore given by
\[
    \Delta_u\longrightarrow\Delta_t,
    \qquad t=u^r.
\]
Consequently, $\widetilde B_O$ has a unique point
above $0$, with ramification index $r$.

When $r=1$, the normalization maps isomorphically
to $\Delta$, and its map to $\mathcal E$ is a section
of $\pi$. Conversely, a section has degree one
over $\Delta$, so its monodromy orbit has
cardinality one.
\end{proof}

Thus monodromy determines the irreducible
components of $T_n$, their degrees over
$\Delta$, and their ramification indices
above $0$. In the following sections,
we determine their specialization points,
their intersections with the components
of the central fiber, and which
irreducible components of $T_n$ meet
at the same point.

\section{Torsion Curves over Fibers of Type \texorpdfstring{$I_k$}{I\_k}}
\label{sec:Ik}

We first consider fibers of type $I_k$, with $k\geq1$.
Their local monodromy is unipotent, and a ramified
base change of degree $h$, followed by minimal
resolution, gives a fiber of type $I_{kh}$.
The monodromy action determines the irreducible
components of $T_n$, while the base change allows
us to compute their specialization points and
intersection multiplicities.

Write
\[
    \mathcal E_0=C_0+\cdots+C_{k-1},
\]
where the components are labeled cyclically and
the zero section meets $C_0$. We write $P_1,\ldots,P_k$ for the nodes in cyclic order, with
$P_j$ lying between the branches labeled $C_{j-1}$ and $C_j$ (indices
modulo $k$). For $k=1$, both branches at $P_1$ belong to $C_0$;
for $k=2$, the two nodes joining $C_0$ and $C_1$ are distinct. As in
the introduction, put
$C_j^*:=C_j\setminus\operatorname{Sing}(\mathcal E_0)$.

When several semistable models occur in the same argument, we write $\mathcal E^{I_r}$, $T_n^{I_r}$, and $C_j^{I_r}$ to indicate that the central fiber has type \(I_r\).

\subsection{Monodromy Orbits and Irreducible Components of \texorpdfstring{$T_n$}{T\_n}}
\label{subsec:Ik-monodromy}

With a suitable basis of $H_1(E_t,\mathbb Z)$, the local monodromy
matrix is
\[
    M_k=\begin{pmatrix}1&k\\0&1\end{pmatrix};
\]
see \cite{Kodaira1963,Miranda1989}. We choose the monodromy basis compatibly with
the component-group labeling after the
degree-$n$ semistable base change, so that
the second coordinate $i$ corresponds to
the component $C_{ik}^{I_{nk}}$.

As explained in Section~2,
this matrix also describes the monodromy on
$E_t[n]\simeq(\mathbb Z/n\mathbb Z)^2$ after reduction modulo $n$.
For $v=(a,i)^T$, we have $M_kv=(a+ik,i)^T$: the second coordinate is
fixed, and the first is translated by $ik$.

Put $d_i:=\gcd(ik,n)$ for $0\leq i\leq n-1$. The subgroup generated
by $ik$ in $\mathbb Z/n\mathbb Z$ is $d_i\mathbb Z/n\mathbb Z$.
Thus the monodromy orbits are
\[
    O_{i,c}:=\{(a,i):a\equiv c\pmod{d_i}\},
    \qquad c\in\mathbb Z/d_i\mathbb Z,
\]
and each has $n/d_i$ elements. By
Proposition~\ref{prop:monodromy-components}, $O_{i,c}$ corresponds to
an irreducible component $B_{i,c}\subset T_n$.

\begin{proposition}[Monodromy orbits for an \(I_k\)-fiber]
\label{prop:Ik-orbits}
Set $g:=\gcd(n,k)$ and $m:=n/g$. Then:
\begin{enumerate}
\item For \(0\leq i\leq n-1\) and
\(c\in\mathbb Z/d_i\mathbb Z\),
\[
\deg(B_{i,c}/\Delta)
=
\frac{n}{d_i}
=
\frac{m}{\gcd(i,m)}.
\]

\item Exactly \(gn\) irreducible components of \(T_n\) are sections.

\item The number of irreducible components of \(T_n\) is
\[
\#\operatorname{Irr}(T_n)
=
\sum_{i=0}^{n-1}d_i
=
\sum_{d\mid m}
g^2d\,\varphi\!\left(\frac{m}{d}\right),
\]
where $\varphi$ is Euler's totient function.
\end{enumerate}
\end{proposition}

\begin{proof}
The $\ell$-th iterate sends $(a,i)^T$ to $(a+\ell ik,i)^T$.
The orbit closes when $\ell ik\equiv0\pmod n$, so its cardinality is
$n/\gcd(ik,n)=n/d_i$. Write $n=gm$ and $k=gk'$, where
$\gcd(k',m)=1$. Then $d_i=g\gcd(ik',m)=g\gcd(i,m)$, therefore
$n/d_i=m/\gcd(i,m)$. This proves~(1).

An orbit is fixed exactly when $n\mid ik$, equivalently
$m\mid i$. There are $g$ such values of $i$ modulo $n$, and for
each of them all $n$ points $(a,i)$ are fixed. Thus $T_n$ has
$gn$ irreducible components that are sections, proving~(2).

For each $i$ there are $d_i$ orbits. Since each residue modulo $m$
occurs $g$ times among $0,\ldots,n-1$, we obtain
\[
\#\operatorname{Irr}(T_n)
=\sum_{i=0}^{n-1}d_i
=g^2\sum_{a=0}^{m-1}\gcd(a,m)
=g^2\sum_{d\mid m}d\,\varphi\!\left(\frac{m}{d}\right).
\]
In the last equality we grouped the residues $a$ according to
$d=\gcd(a,m)$. This proves~(3).
\end{proof}

In particular, if  $\gcd(n,k)=1$, there are exactly $n$ section components;
the degree of $B_{i,c}$ over $\Delta$ is $n/\gcd(i,n)$.
If $n\mid k$, monodromy acts trivially on $E_t[n]$, and \(T_n\) is a union of \(n^2\) sections.

\subsection{Smooth and Nodal Specialization}
\label{subsec:Ik-smooth-nodal}

The specialization of an irreducible component
of $T_n$ is constrained by its ramification
over $0$. For an $I_k$-fiber, we show that
the section components specialize to smooth
points of the central fiber, while all other
irreducible components specialize to nodes.

\begin{lemma}[Smooth and nodal specialization]
\label{lem:Ik-smooth-nodal}
An irreducible component of \(T_n\) meets the smooth locus of
\(\mathcal E_0\) if and only if its monodromy orbit is fixed. Consequently,
the section components specialize to the loci \(C_j^*\), whereas every
component corresponding to a nontrivial monodromy orbit specializes to one
of the nodes \(P_j\).
\end{lemma}

\begin{proof}
At a node of $\mathcal E_0$ the family has the local equation $xy=t$.
If a section passed through the node, both $x(t)$ and $y(t)$ would
vanish at $t=0$, giving
$\operatorname{ord}_0x+\operatorname{ord}_0y\geq2$, contrary to
$xy=t$. Hence a section meets $\mathcal E_0$ at a smooth point.

Conversely, near a smooth point of \(\mathcal E_0\), the morphism
\(\pi\) is smooth. The smooth locus of the minimal regular model is
the N\'eron model of the generic fiber.
Since multiplication by $n$ on this smooth group
scheme is \'etale in characteristic zero,
the torsion locus is \'etale over the base
near such a point. Hence any irreducible component specializing there is
unramified over \(0\). By Proposition~\ref{prop:monodromy-components},
its monodromy orbit has length one.
\end{proof}


For a fiber of type $I_k$, the identity component
of the N\'eron model is $\mathbb G_m$, and its
component group is $\Phi \simeq \mathbb Z/k\mathbb Z$;
see \cite[Chap.~IV]{SilvermanAdvanced}. Thus
\[
\mathcal E_0^{\mathrm{sm}}
=
\bigsqcup_{j=0}^{k-1} C_j^*,
\]
where each \(C_j^*\) is a torsor under \(\mathbb G_m\).

\begin{lemma}
\label{lem:Ik-smooth-bound}
For every \(0\leq j\leq k-1\),
\[
\operatorname{length}(T_n\cap C_j^*)\leq n.
\]
\end{lemma}

\begin{proof}
If \(C_j^*\) contains a specialization of an \(n\)-torsion point, then the
\(n\)-torsion lying on \(C_j^*\) is a torsor under the \(n\)-torsion subgroup
of the identity component. Since
\[
\mathbb G_m[n]=\mu_n
\]
has length \(n\) in characteristic zero, the intersection
\(T_n\cap C_j^*\) has length at most \(n\).
\end{proof}

\subsection{Ramified Base Change and Exact Specialization}
\label{subsec:base-change-Ik}
We use ramified base change to determine
the specialization point of each $B_{i,c}$.
For $h\geq1$, consider the base change
$t=s^h$. Locally at a node of an $I_k$-fiber,
the equation $xy=t$ becomes $xy=s^h$. For $h>1$ this is an $A_{h-1}$ surface singularity, whose
minimal resolution inserts a chain of $h-1$ rational curves.
Resolving the singularities above the $k$ nodes gives a central
fiber of type $I_{kh}$ \cite[Lecture VI, \S 4, Table VI.4.1]{Miranda1989}.

\begin{proposition}[Ramified base change of an \(I_k\)-fiber]
\label{prop:base-change-Ik}
Let $f\colon\mathcal E^{I_{kh}}\to\mathcal E^{I_k}$ be the map
obtained by the degree-$h$ base change, followed by normalization
and minimal resolution. Label the components of the new central
fiber cyclically, so that $C_{ih}^{I_{kh}}$ is the strict transform
of $C_i^{I_k}$. Then
\[
 f(C_{ih}^{I_{kh}})=C_i^{I_k},\qquad
 f(C_{ih+r}^{I_{kh}})=P_{i+1}\quad(1\leq r<h).
\]
\end{proposition}

\begin{proof}
The resolution of $xy=s^h$ inserts $h-1$ exceptional curves
between the two local branches over each node (no exceptional curve
when $h=1$). In the cyclic labeling, these are
$C_{ih+1}^{I_{kh}},\ldots,C_{ih+h-1}^{I_{kh}}$ above $P_{i+1}$.
They are contracted by $f$, while $C_{ih}^{I_{kh}}$ maps onto
$C_i^{I_k}$, as claimed.
\end{proof}

To determine the specialization of \(T_n\), we first consider two cases in which the torsion curve splits completely into sections. 

Suppose first that the central fiber is of type $I_n$. The monodromy acts trivially on $n$-torsion,
so $T_n$ consists of $n^2$ sections.
By Lemma~\ref{lem:Ik-smooth-bound}, at most
$n$ sections specialize to the smooth locus
of each component of the central fiber.
Since all $n^2$ sections specialize to
smooth points and the fiber has $n$
irreducible components, exactly $n$
sections meet each component.

Next consider a fiber of type $I_{nk}$,
obtained from an $I_n$-fiber by a ramified
base change of degree $k$.
By Proposition~\ref{prop:base-change-Ik},
the components $C_{jk}^{I_{nk}}$,
for $0\leq j\leq n-1$, are the strict
transforms of the components of the
$I_n$-fiber. The remaining components
are contracted to its nodes.
Consequently, the $n$-torsion sections
meet precisely the components
$C_{jk}^{I_{nk}}$, with $n$ sections
meeting each of them.

Returning to a fiber of type $I_k$, apply
the degree-$n$ base change $t=s^n$,
followed by normalization and minimal
resolution:
\[
    f:\mathcal E^{I_{nk}}
      \longrightarrow\mathcal E^{I_k}.
\]
The monodromy on $n$-torsion becomes
trivial, so the pulled-back torsion cover
splits into $n^2$ sections. The preceding discussion
determines which components of the \(I_{nk}\)-fiber these sections meet,
and Proposition~\ref{prop:base-change-Ik} determines their images on the
original \(I_k\)-fiber.

\begin{proposition}[Specialization of an individual component]
\label{prop:Ik-branch-specialization}
Let \(B_{i,c}\subset T_n^{I_k}\) correspond to the monodromy orbit
\[
O_{i,c}
=
\{(a,i):a\equiv c\pmod{d_i}\},
\qquad
d_i=\gcd(ik,n).
\]
After the degree-$n$ base change $t=s^n$ and resolution, the
pullback of $B_{i,c}$ splits into torsion sections, each meeting
$(C_{ik}^{I_{nk}})^*$. Here $f\colon\mathcal E^{I_{nk}}\to\mathcal E^{I_k}$
denotes the map of Proposition~\ref{prop:base-change-Ik}.

If \(n\mid ik\), write \(ik=jn\) with \(0\leq j\leq k-1\). Then
\(B_{i,c}\) is a section and specializes to a point of \(C_j^*\).

If $n\nmid ik$, let $j\in\{1,\ldots,k\}$ be the unique integer such that
$(j-1)n<ik<jn$. Then $B_{i,c}$ specializes to $P_j$.
\end{proposition}

\begin{proof}
Suppose first that \(ik=jn\). Then
\(C_{ik}^{I_{nk}}=C_{jn}^{I_{nk}}\) is the strict transform of
\(C_j^{I_k}\). Moreover, \(d_i=\gcd(ik,n)=n\), so
\(O_{i,c}\) has length \(n/d_i=1\). Hence \(B_{i,c}\) is a
section, and by Lemma~\ref{lem:Ik-smooth-nodal} it specializes to a point
of \(C_j^*\).

If $(j-1)n<ik<jn$, the component $C_{ik}^{I_{nk}}$ lies
in the exceptional chain above $P_j$ and is contracted to that node
by Proposition~\ref{prop:base-change-Ik}. Every lifted torsion
section meets this component, so $B_{i,c}$ specializes to $P_j$.
\end{proof}

We now sum the contributions of the components
$B_{i,c}$ to obtain the intersection lengths
at the nodes and in the smooth locus of
$\mathcal E_0$.

\begin{theorem}[Specialization over an \(I_k\)-fiber]
\label{thm:torsion-distribution-Ik}
Let \(\mathcal E_0\) be of type \(I_k\). Then:

\begin{enumerate}
\item For \(0\leq j\leq k-1\),
\[
\operatorname{length}(T_n\cap C_j^*)
=
\begin{cases}
n,& k\mid jn,\\
0,& k\nmid jn.
\end{cases}
\]

\item For \(1\leq j\leq k\),
\[
\operatorname{length}_{P_j}(T_n\cap\mathcal E_0)
=
n\, \cdot
\#  \left\{
1\leq i\leq n-1:
(j-1)n<ik<jn
\right\}.
\]

\item For every \(i\) satisfying \((j-1)n<ik<jn\), set
\(d_i=\gcd(ik,n)\). There are $d_i$ irreducible components $B_{i,c}$, indexed
by $c\in\mathbb Z/d_i\mathbb Z$, specializing to $P_j$.
Each has a unique analytic branch there, with ramification index
$n/d_i$ over $0$.
\end{enumerate}
\end{theorem}

\begin{proof}
By Proposition~\ref{prop:Ik-branch-specialization}, a component specializes
to \(C_j^*\) precisely when \(ik=jn\) for some \(0\leq i\leq n-1\).
Such an \(i\) exists if and only if \(k\mid jn\). In this case
\(d_i=\gcd(ik,n)=n\), so there are \(n\) singleton monodromy orbits with
second coordinate \(i\). Hence exactly \(n\) section components specialize
to \(C_j^*\), proving~(1).

Fix a node \(P_j\). By
Proposition~\ref{prop:Ik-branch-specialization}, the components specializing
to \(P_j\) are precisely those indexed by \(i\) satisfying
\((j-1)n<ik<jn\). For each such \(i\), there are
\(d_i=\gcd(ik,n)\) components \(B_{i,c}\), each of ramification index
\(n/d_i\). This proves~(3).

The intersection multiplicity of each such component with
$\mathcal E_0$ equals the degree $n/d_i$ of its normalization over
$\Delta$. Hence the $d_i$ components indexed by $i$ contribute
$d_i(n/d_i)=n$ to the local intersection length. Summing over the
admissible values of $i$ gives
\[
\operatorname{length}_{P_j}(T_n\cap\mathcal E_0)
=
n\,
\#\left\{
1\leq i\leq n-1:
(j-1)n<ik<jn
\right\},
\]
proving~(2).
\end{proof}

This theorem also determines the connected components of the
closed curve \(T_n\). The section components are pairwise disjoint and
specialize to smooth points of the central fiber. By
Proposition~\ref{prop:Ik-orbits}, there are \(gn\) such components, where
\(g=\gcd(n,k)\).

Every remaining irreducible component specializes to a unique node \(P_j\),
and all components specializing to the same node meet there. Thus each node
met by the non-section part of \(T_n\) contributes one connected component.
By Theorem~\ref{thm:torsion-distribution-Ik}, \(P_j\) is met precisely when
the interval \(((j-1)n,jn)\) contains a multiple of \(k\).

\begin{corollary}[Connected components of \(T_n\)]
\label{cor:connected-components-Tn}
Let \(g=\gcd(n,k)\). Then
\[
\#\pi_0(T_n)
=
gn+\min\{n-g,k\}.
\]
\end{corollary}

\begin{proof}
It remains to count the nodes met by the non-section components. If \(n>k\), every interval \(((j-1)n,jn)\) has length greater than \(k\),
and hence contains a multiple of \(k\). Thus all \(k\) nodes are met.

Suppose \(n\leq k\). Each interval has length at most \(k\), so it contains
at most one multiple of \(k\). The relevant multiples are
\(k,2k,\ldots,(n-1)k\). A multiple \(ik\) lies on an endpoint of one of
the intervals precisely when \(n\mid ik\). Since \(g=\gcd(n,k)\), this is
equivalent to \((n/g)\mid i\). Among \(1\leq i\leq n-1\), there are exactly
\(g-1\) such values. Therefore
\(
(n-1)-(g-1)=n-g
\)
multiples lie in the interiors of the intervals, so exactly \(n-g\) nodes
are met.

Thus the number of nodes met is \(\min\{n-g,k\}\). Adding the \(gn\)
section components gives
\[
\#\pi_0(T_n)
=
gn+\min\{n-g,k\}.
\]
\end{proof}

\subsection{Intersection Numbers of Irreducible Components}
\label{subsec:branchwise-intersections-Ik}

The degree-$n$ base change also allows us
to compute the intersection numbers
$B_{i,c}\cdot C_j$. We first determine the
pullbacks of the divisors $C_j$ to the
resolved $I_{nk}$-model.

\begin{lemma}[Pullback of the fiber components]
\label{lem:pullback-components-Ik}
For every \(j\),
\[
f^*(C_j^{I_k})
=
nC_{nj}^{I_{nk}}
+
\sum_{r=1}^{n-1}(n-r)
\left(
C_{nj-r}^{I_{nk}}+C_{nj+r}^{I_{nk}}
\right),
\]
where the indices on the \(I_{nk}\)-fiber are taken modulo \(nk\).
\end{lemma}

\begin{proof}
First suppose $k \geq 2$. We compute the total transform of \(C_j^{I_k}\) near its two nodes. Consider first
\(
P_{j+1}.
\)
Since the total space is smooth and the central fiber has simple normal
crossings, there are étale-local coordinates \(x,y\) and a parameter \(t\)
on the base such that
\[
xy=t,
\qquad
C_j^{I_k}=\{x=0\},
\qquad
C_{j+1}^{I_k}=\{y=0\}.
\]
After the degree-\(n\) base change \(t=s^n\), the local equation becomes
\(xy=s^n\), an \(A_{n-1}\)-surface singularity.

We describe its minimal resolution torically. The surface \(xy=s^n\)
corresponds to the cone generated by
\[
v_0=(1,0),
\qquad
v_n=(1,n),
\]
and its minimal regular subdivision is obtained by inserting the rays
\(v_r=(1,r)\), \(1\leq r\leq n-1\); see
\cite[Chaps.~1 and~3]{CoxLittleSchenck2011}.
The corresponding divisors form the exceptional chain.

The function \(x\) is the torus character corresponding to
\(m_x=(n,-1)\). For a torus character \(\chi^m\), the order of vanishing
along the torus-invariant divisor associated with a primitive ray \(v\) is
\[
\operatorname{ord}_{D_v}(\chi^m)=\langle m,v\rangle;
\]
see \cite[Chap.~4]{CoxLittleSchenck2011}. Hence
\[
\operatorname{ord}_{D_{v_r}}(x)
=
\langle(n,-1),(1,r)\rangle
=
n-r.
\]

The total transform of \(C_j^{I_k}=\{x=0\}\) near \(P_{j+1}\) is
\[
nC_{nj}^{I_{nk}}
+
\sum_{r=1}^{n-1}(n-r)C_{nj+r}^{I_{nk}}.
\]

At the other node \(P_j=C_{j-1}^{I_k}\cap C_j^{I_k}\), the same calculation,
with the exceptional chain reversed, gives
\[
nC_{nj}^{I_{nk}}
+
\sum_{r=1}^{n-1}(n-r)C_{nj-r}^{I_{nk}}.
\]
Combining the two contributions and counting the strict transform
\(C_{nj}^{I_{nk}}\) only once gives the stated formula.

For $k=1$, the central fiber is the single nodal divisor
$C_0^{I_1}=\mathcal E_0^{I_1}$. Thus
\[f^*C_0^{I_1}=n\mathcal E_0^{I_n}=n\sum_{a=0}^{n-1}C_a^{I_n}.\]
In the formula of the lemma, each exceptional component occurs
twice, with coefficients $r$ and $n-r$, so the same identity holds.
\end{proof}

The pullback formula determines the intersection multiplicities of each
irreducible component \(B_{i,c}\subset T_n\) with the components of the
central fiber.

\begin{theorem}[Intersection formula for irreducible components]
\label{thm:branchwise-intersections-Ik}
Let \(B_{i,c}\subset T_n^{I_k}\), and set \(d_i=\gcd(ik,n)\).

Suppose \(k\geq2\) and
\(
(j-1)n<ik<jn.
\)
Then \(B_{i,c}\) specializes to \(P_j=C_{j-1}\cap C_j\), and
\[
B_{i,c}\cdot C_{j-1}
=
\frac{jn-ik}{d_i},
\qquad
B_{i,c}\cdot C_j
=
\frac{ik-(j-1)n}{d_i}.
\]
Moreover, \(B_{i,c}\cdot C_r=0\) for every \(r\neq j-1,j\).

If \(ik=jn\), then \(d_i=n\), the component \(B_{i,c}\) is a section
specializing to \(C_j^*\), and
\[
B_{i,c}\cdot C_j=1,
\qquad
B_{i,c}\cdot C_r=0
\quad (r\neq j).
\]
For $k=1$ and $1\leq i\leq n-1$, the nodal specialization instead
satisfies $B_{i,c}\cdot C_0=n/d_i$.
\end{theorem}

\begin{proof}
Suppose first that $k \geq 2$ and \((j-1)n<ik<jn\). After the degree-\(n\) base change,
\(B_{i,c}\) splits into torsion sections, each meeting \(C_{ik}^{I_{nk}}\). Choose one such section and denote it by \(S_{a,i}\).
The induced map \(S_{a,i}\to B_{i,c}\) has degree \(d_i\), so
\[
f_*S_{a,i}=d_iB_{i,c}.
\]
By the projection formula,
\[
S_{a,i}\cdot f^*C_{j-1}
=
d_i(B_{i,c}\cdot C_{j-1}),
\qquad
S_{a,i}\cdot f^*C_j
=
d_i(B_{i,c}\cdot C_j).
\]

The section \(S_{a,i}\) meets \(C_{ik}^{I_{nk}}\) transversely at a smooth point of the central fiber. Hence its intersection with a pullback divisor is the coefficient of \(C_{ik}^{I_{nk}}\) in that divisor. By Lemma~\ref{lem:pullback-components-Ik}, these coefficients are \(jn-ik\) in \(f^*C_{j-1}\) and \(ik-(j-1)n\) in \(f^*C_j\). Therefore
\[
B_{i,c}\cdot C_{j-1}
=
\frac{jn-ik}{d_i},
\qquad
B_{i,c}\cdot C_j
=
\frac{ik-(j-1)n}{d_i}.
\]
The component \(C_{ik}^{I_{nk}}\) does not occur in the pullback of any
other \(C_r\), so \(B_{i,c}\cdot C_r=0\) for \(r\neq j-1,j\).

Now suppose \(ik=jn\). Then \(d_i=n\), and
\(C_{ik}^{I_{nk}}=C_{jn}^{I_{nk}}\) is the strict transform of
\(C_j^{I_k}\). By Lemma~\ref{lem:pullback-components-Ik}, it occurs with
coefficient \(n\) in \(f^*C_j\). The projection formula therefore gives
\[
B_{i,c}\cdot C_j
=
\frac{n}{d_i}
=
1.
\]
Since \(B_{i,c}\) is a section specializing to \(C_j^*\), it is disjoint
from every other component of the central fiber.

Finally, if $k=1$, we have $C_0=\mathcal E_0$ as a divisor.
Hence $B_{i,c}\cdot C_0=\deg(B_{i,c}/\Delta)=n/d_i$.
\end{proof}

\subsection{Component-wise Intersection Numbers}
\label{subsec:intersections-Ik}

The intersection numbers of the full torsion curve with the components of
the \(I_k\)-fiber are obtained by summing the contributions of the individual
components \(B_{i,c}\). Put
$(x)_+:=\max\{x,0\}$.

\begin{theorem}[Component-wise intersection formula]
\label{thm:intersection-numbers-Ik}
For every \(0\leq j\leq k-1\),
\[
T_n\cdot C_j
=
\sum_{r\in\mathbb Z}
\bigl(n-|nj-rk|\bigr)_+.
\]
Only finitely many terms in the sum are nonzero.
\end{theorem}

\begin{proof}
Suppose first that $k\geq2$. Fix $i\in\{0,\ldots,n-1\}$.
For $(j-1)n<ik<jn$, the $d_i$ components indexed by $i$
contribute $ik-(j-1)n$ to $T_n\cdot C_j$. For
$jn<ik<(j+1)n$, their total contribution is $(j+1)n-ik$.
When $ik=jn$, the $n$ corresponding sections contribute $n$.
Thus the total contribution for this $i$ is
$(n-|nj-ik|)_+$. Using the cyclic indexing of the $I_k$-fiber,
we obtain
\[
    T_n\cdot C_j=\sum_{i\in\mathbb Z}(n-|nj-ik|)_+.
\]
For $k=1$, the central fiber is $C_0$ and
$T_n\cdot C_0=n^2$; the same formula holds because
$\sum_{i\in\mathbb Z}(n-|i|)_+=n^2$.
\end{proof}

Summing over the fiber components gives
\[\sum_{j=0}^{k-1}T_n\cdot C_j=T_n\cdot\mathcal E_0=n^2,\]
as expected from the degree of $T_n\to\Delta$.

The next formula shows that replacing $n$ by $n+k$ adds the
same amount to the intersection number with every $C_j$.

\begin{theorem}[Reduction modulo \(k\)]
\label{thm:intersection-numbers-Ik-periodicity}
Write $n=\alpha k+\beta$, with $0\leq\beta<k$. Then,
for $0\leq j<k$,
\[
 T_n\cdot C_j
 = T_{\beta}\cdot C_j
  +\frac{n^2-\beta^2}{k}.
\]

Here we set $T_{0}\cdot C_j=0$.
\end{theorem}

\begin{proof}
For $r\geq0$, put
\[
 S_r(x):=\sum_{q\in\mathbb Z}(r-|x-qk|)_+.
\]
Theorem~\ref{thm:intersection-numbers-Ik} gives
$T_n\cdot C_j=S_n(nj)=S_n(\beta j)$, since $S_n$ is
$k$-periodic.

For $r\geq0$, the difference $S_{r+k}-S_r$ is constant:
away from the break points, the two piecewise-linear functions
have the same derivative. Moreover,
\[
 \int_0^k S_r(x)\,dx
 =\int_{\mathbb R}(r-|u|)_+\,du=r^2.
\]
Integration over one period therefore gives
\[S_{r+k}(x)-S_r(x)=((r+k)^2-r^2)/k.\]
Iterating from $r=\beta$ to $r=n$ yields
\[
 S_n(x)-S_\beta(x)=\frac{n^2-\beta^2}{k}.
\]
Set $x=\beta j$ and expand $S_\beta(\beta j)$ to obtain
the claimed formula.
\end{proof}

Combining Proposition~\ref{prop:Ik-orbits},
Theorem~\ref{thm:torsion-distribution-Ik},
Corollary~\ref{cor:connected-components-Tn},
Theorem~\ref{thm:branchwise-intersections-Ik}, and
Theorem~\ref{thm:intersection-numbers-Ik-periodicity}
completes the analysis of the \(I_k\)-case.

\section{Torsion Curves over Fibers of Type
\texorpdfstring{$I_k^*$}{I\_k*}}
\label{sec:torsion-Ik-star}

We now consider fibers of type $I_k^*$, with $k\geq1$.
Their local monodromy is $-M_k$, so the quadratic
base change $t=s^2$ gives monodromy $M_{2k}$ and,
after resolution, a fiber of type $I_{2k}$.
The change $s\mapsto-s$ acts on this new model,
and its action on the components of the central
fiber allows us to recover the $I_k^*$-configuration
and determine the specialization of the torsion curve.

Write \(D_0,\ldots,D_k\) for the multiplicity-two components of the
central chain, and put \(P_j:=D_{j-1}\cap D_j\) for \(1\leq j\leq k\).
The four multiplicity-one components are \(D_a,D_b,D_c,D_d\), where
\(D_a,D_b\) meet \(D_0\) and \(D_c,D_d\) meet \(D_k\). We label them so
that the zero section meets \(D_a\). For every irreducible component
\(D\subset\mathcal E_0\), put
\(D^*:=D\setminus\operatorname{Sing}((\mathcal E_0)_{\mathrm{red}})\).
Thus \(D^*\) excludes the nodes of the reduced central fiber; it is
not the smooth locus of the fibration when \(D\) has multiplicity two.

The outer components have multiplicity one, and the components
\(D_0,\ldots,D_k\) have multiplicity two. If an irreducible component
\(B\subset T_n\) of degree \(r\) over \(\Delta\) specializes to a point
of \(D^*\), then
\[
r
=
B\cdot\mathcal E_0
=
m_D(B\cdot D),
\]
where \(m_D\) denotes the multiplicity of \(D\) in \(E_0\).
When both models occur in the same argument, we use
\(\mathcal E^{I_k^*},T_n^{I_k^*}\) for the original model and
\(\mathcal E^{I_{2k}},T_n^{I_{2k}}\) for the semistable model after
quadratic base change.

For reference, we collect the results proved below in a single theorem.
For any condition \(\mathcal P\), let \(\mathds{1}_{\mathcal P}\) denote
its indicator.

\begin{theorem}[Torsion over an \(I_k^*\)-fiber]
\label{thm:Ik-star-complete}
Let \(\mathcal E_0\) be of type \(I_k^*\), \(k\geq1\), with the
component labeling above, and set \(g:=\gcd(n,2k)\). Then:

\begin{enumerate}
\item The number of irreducible components of \(T_n\) is
\[
\#\operatorname{Irr}(T_n)
=
\frac12
\sum_{m\mid n}
\varphi(m)\frac{n}{m}\gcd(m,2k)
+
\frac{\gcd(n,2)}{2}
+
\mathds{1}_{\,2\mid n}\mathds{1}_{\,4\mid nk}.
\]

The number of section components is
\[
\begin{cases}
1,&n\text{ odd},\\
2,&2\mid n\text{ and }4\nmid nk,\\
4,&4\mid nk.
\end{cases}
\]

\item The intersection lengths on the outer components away from the nodes are
\[
\operatorname{length}(T_n\cap D_a^*)=1,
\qquad
\operatorname{length}(T_n\cap D_b^*)=
\begin{cases}
1,&2\mid n,\\
0,&2\nmid n,
\end{cases}
\]
and
\[
\operatorname{length}(T_n\cap D_c^*)
=
\operatorname{length}(T_n\cap D_d^*)
=
\begin{cases}
1,&2\mid n\text{ and }4\mid nk,\\
0,&\text{otherwise}.
\end{cases}
\]
For the multiplicity-two components,
\[
\operatorname{length}(T_n\cap D_0^*)
=
\begin{cases}
\dfrac{n-1}{2},&2\nmid n,\\[0.6em]
\dfrac{n-2}{2},&2\mid n,
\end{cases}
\]
\[
\operatorname{length}(T_n\cap D_k^*)
=
\begin{cases}
\dfrac{n-2}{2},&2\mid n\text{ and }4\mid nk,\\[0.6em]
\dfrac n2,&2\mid n\text{ and }4\nmid nk,\\[0.6em]
0,&2\nmid n,
\end{cases}
\]
and, for \(1\leq j\leq k-1\),
\[
\operatorname{length}(T_n\cap D_j^*)
=
\begin{cases}
n,&2k\mid jn,\\
0,&2k\nmid jn.
\end{cases}
\]

\item The torsion curve avoids the four outer nodes. For
\(1\leq j\leq k\), the node \(P_j\) is met precisely for those
\(i\) satisfying
\[
(j-1)n<2ki<jn.
\]
For every such \(i\), put \(d_i:=\gcd(2ki,n)\). There are
\(d_i\) irreducible components \(B_{i,\ell}\), each with a unique analytic
branch at \(P_j\) and ramification index \(2n/d_i\) over \(0\). They satisfy
\[
B_{i,\ell}\cdot D_{j-1}
=
\frac{jn-2ki}{d_i},
\qquad
B_{i,\ell}\cdot D_j
=
\frac{2ki-(j-1)n}{d_i}.
\]
In particular,
\[
\operatorname{length}_{P_j}(T_n\cap\mathcal E_0)
=
2n\,
\#\left\{
1\leq i\leq n-1:
(j-1)n<2ki<jn
\right\}.
\]

\item For \(1\leq r\leq k-1\),
\[
T_n\cdot D_r
=
\sum_{i\in\mathbb Z}
\bigl(n-|nr-2ki|\bigr)_+.
\]
At the two ends of the multiplicity-two chain,
\[
T_n\cdot D_0
=
\frac12
\sum_{i\in\mathbb Z}
\bigl(n-|2ki|\bigr)_+
-
\frac{\gcd(n,2)}{2},
\]
and
\[
T_n\cdot D_k
=
\frac12
\sum_{i\in\mathbb Z}
\bigl(n-|nk-2ki|\bigr)_+
-
\mathds{1}_{\,2\mid n}\mathds{1}_{\,4\mid nk}.
\]
For \(\alpha\in\{a,b,c,d\}\),
\[
T_n\cdot D_\alpha= \operatorname{length}(T_n\cap D_\alpha^*).
\]

\item The number of connected components is
\[
\#\pi_0(T_n)
=
\frac{ng+\gcd(n,2)}{2}
+
\mathds{1}_{\,2\mid n}\mathds{1}_{\,4\mid nk}
+
\min\left\{k,\frac{n-g}{2}\right\}.
\]
\end{enumerate}
\end{theorem}

\begin{proof}
Part~(1) follows from
Propositions~\ref{prop:section-components-Ik-star} and
\ref{prop:irreducible-components-Ik-star}.
Part~(2) is Theorem~\ref{thm:smooth-specialization-Ik-star}.
Part~(3) follows from
Proposition~\ref{prop:outer-components-Ik-star} and
Theorem~\ref{thm:branches-Ik-star}.
Part~(4) is Corollary~\ref{cor:componentwise-intersections-Ik-star},
and part~(5) is
Proposition~\ref{prop:connected-components-Ik-star}.
\end{proof}

\subsection{Monodromy Orbits}
\label{subsec:Ik-star-monodromy}

With respect to a suitable basis of \(H_1(E_t,\mathbb Z)\), the
local monodromy of an \(I_k^*\)-fiber is
\[
-M_k
=
-\begin{pmatrix}
1&k\\
0&1
\end{pmatrix};
\]
see \cite{Kodaira1963,Miranda1989}. Hence its action on
\(E_t[n]\simeq(\mathbb Z/n\mathbb Z)^2\) is
\(
(x,y)\longmapsto(-x-ky,-y).
\)

We first determine the fixed points. They satisfy
\[
2y=0,
\qquad
2x+ky=0
\quad\text{in }\mathbb Z/n\mathbb Z.
\]
If \(n\) is odd, the only solution is \((0,0)\). If \(n\) is even, the
solutions with \(y=0\) are the two points satisfying \(2x=0\). For
\(y=n/2\), solutions exist if and only if \(4\mid nk\), and in that case
there are again two.

\begin{proposition}[Section components]
\label{prop:section-components-Ik-star}
The number of section components of \(T_n\) is
\[
\begin{cases}
1,& n\text{ odd},\\
2,& 2\mid n\text{ and }4\nmid nk,\\
4,& 4\mid nk.
\end{cases}
\]
\end{proposition}

\begin{proof}
By Proposition~\ref{prop:monodromy-components}, section components
correspond exactly to fixed monodromy orbits, so the result follows from the
preceding calculation.
\end{proof}

To count all monodromy orbits, observe that
\[
(-M_k)^2
=
\begin{pmatrix}
1&2k\\
0&1
\end{pmatrix},
\qquad
(-M_k)^2(x,y)=(x+2ky,y).
\]
Thus the square of the monodromy preserves the second coordinate and
translates the first coordinate by \(2ky\).

\begin{proposition}[Number of irreducible components]
\label{prop:irreducible-components-Ik-star}
The number of irreducible components of \(T_n\) is
\[
\#\operatorname{Irr}(T_n)
=
\frac12
\sum_{m\mid n}
\varphi(m)\frac{n}{m}\gcd(m,2k)
+
\frac{\gcd(n,2)}{2}
+
\mathds{1}_{\,2\mid n}\mathds{1}_{\,4\mid nk}.
\]
\end{proposition}

\begin{proof}
Let
\[
A_y
:=
\{(x,y):x\in\mathbb Z/n\mathbb Z\}
\subset E_t[n].
\]
The monodromy sends \(A_y\) to \(A_{-y}\), while its square preserves
\(A_y\) and acts there by
\(
x\longmapsto x+2ky.
\)

Suppose that \(y\) has exact order \(m\). Writing
\(y=(n/m)u\), with \(u\in(\mathbb Z/m\mathbb Z)^\times\), the translation
by \(2ky\) has
\[
\gcd(n,2ky)
=
\frac{n}{m}\gcd(m,2k)
\]
orbits on \(A_y\).

We now distinguish whether \(y=-y\).

If \(m>2\), then \(y\neq-y\). The monodromy exchanges \(A_y\) and
\(A_{-y}\), so each orbit of its square on \(A_y\) combines with the
corresponding orbit on \(A_{-y}\) to form one monodromy orbit on
\(A_y\cup A_{-y}\). Thus the pair \(\{y,-y\}\) contributes
\[
\frac{n}{m}\gcd(m,2k)
\]
monodromy orbits. Since the \(\varphi(m)\) elements of exact order \(m\)
form \(\varphi(m)/2\) pairs \(\{y,-y\}\), their total contribution is
\[
\frac12\varphi(m)\frac{n}{m}\gcd(m,2k).
\]

It remains to consider the self-opposite values of \(y\). The first is
\(y=0\). On \(A_0\), the monodromy is the involution \(x\mapsto-x\).
It has \(\gcd(n,2)\) fixed points, and therefore
\[
\frac{n+\gcd(n,2)}{2}
\]
orbits.

If \(n\) is even, the second self-opposite value is \(y=n/2\). On
\(A_{n/2}\), the monodromy is
\(
x\longmapsto-x-k\frac n2.
\)
Its fixed points satisfy
\[
2x\equiv-k\frac n2\pmod n.
\]
There are two such points precisely when \(4\mid nk\), and none otherwise.
Hence \(A_{n/2}\) contributes
\[
\frac n2+\mathds{1}_{4\mid nk}
\]
monodromy orbits.

The divisor sum
\[
\frac12
\sum_{m\mid n}
\varphi(m)\frac{n}{m}\gcd(m,2k)
\]
already contains \(n/2\) for \(m=1\), and, when \(2\mid n\), another
\(n/2\) for \(m=2\). Replacing these two terms by the actual contributions
of the self-opposite values gives
\[
\#\operatorname{Irr}(T_n)
=
\frac12
\sum_{m\mid n}
\varphi(m)\frac{n}{m}\gcd(m,2k)
+
\frac{\gcd(n,2)}{2}
+
\mathds{1}_{\,2\mid n}\mathds{1}_{\,4\mid nk}.
\]
\end{proof}

\subsection{Quadratic Base Change}
\label{subsec:quadratic-reduction-Ik-star}

Consider the quadratic base change \(q:\Delta_s\to\Delta_t\),
\(t=s^2\). Pulling back \(\pi:\mathcal E^{I_k^*}\to\Delta_t\) and
passing to the relatively minimal regular model gives the semistable
family
\[
\pi'\colon\mathcal E^{I_{2k}}\longrightarrow\Delta_s.
\]

A positive generator of \(\pi_1(\Delta_s^*)\) maps to twice a positive
generator of \(\pi_1(\Delta_t^*)\). Hence the monodromy of the pulled-back
family is
\[
(-M_k)^2=M_{2k}
=
\begin{pmatrix}
1&2k\\
0&1
\end{pmatrix}.
\]
Thus the central fiber of \(\mathcal E^{I_{2k}}\) is of type \(I_{2k}\);
see \cite[Lecture VI, \S 4, Table VI.4.1]{Miranda1989}.

Since $t=s^2$ is unchanged when $s$ is replaced
by $-s$, this change of coordinate induces
an automorphism
\[
    \tau:\mathcal E^{I_{2k}}
    \longrightarrow\mathcal E^{I_{2k}},
    \qquad \tau^2=\mathrm{id}.
\]
We describe how $\tau$ acts on the components
of the central fiber.

\begin{proposition}[Quadratic semistable reduction]
\label{prop:quadratic-reduction-Ik-star}
After a suitable cyclic labeling of the components
\[
C_0^{I_{2k}},\ldots,C_{2k-1}^{I_{2k}},
\]
the following hold:
\begin{enumerate}
\item The involution  $\tau$ acts by
\[
\tau(C_j^{I_{2k}})
=
C_{2k-j}^{I_{2k}},
\]
with indices taken modulo \(2k\).

\item The two opposite components \(C_0^{I_{2k}}\) and
\(C_k^{I_{2k}}\) are preserved by \(\tau\).

\item For $1\leq j\leq k-1$, the map $\tau$
exchanges $C_j^{I_{2k}}$ and
$C_{2k-j}^{I_{2k}}$. These two components
have the same image in the quotient,
corresponding to $D_j$.

\item The images of \(C_0\) and \(C_k\) correspond to \(D_0\) and \(D_k\), respectively. The quotient has singularities
at the images of the fixed points of $\tau$.
Resolving these singularities produces the
four outer components $D_a,D_b,D_c,D_d$.
\end{enumerate}
\end{proposition}

\begin{proof}
The involution $\tau$ preserves the zero section and hence
the identity component $C_0^{I_{2k}}$.
On its smooth locus, identified with $\mathbb G_m$,
the involution acts by $z\mapsto z^{-1}$.
This is the involution arising from the quadratic twist
description of an $I_k^*$-fiber
\cite[Lecture VI, \S 4]{Miranda1989}.
It exchanges the two nodes of $C_0^{I_{2k}}$ and therefore
reverses the cyclic order of the components.
With our labeling,
\[
\tau(C_j^{I_{2k}})=C_{2k-j}^{I_{2k}},
\]
where the indices are taken modulo $2k$.
This proves (1) and (2). For $1\le j\le k-1$, the components
$C_j^{I_{2k}}$ and $C_{2k-j}^{I_{2k}}$
are exchanged by $\tau$ and hence have the
same image in the quotient, giving $D_j$.
This proves (3).

Each preserved component contains two fixed points of \(\tau\).
Their images are quotient singularities, whose resolutions introduce
two multiplicity-one components at each end of the chain. The images of
\(C_0^{I_{2k}}\) and \(C_k^{I_{2k}}\) become \(D_0\) and \(D_k\),
respectively. This proves~(4).
\end{proof}

\subsection{Specialization to the Smooth Locus}
\label{subsec:smooth-specialization-Ik-star}

The four outer nodes of the \(I_k^*\)-fiber are
\[
Q_a=D_a\cap D_0,\qquad
Q_b=D_b\cap D_0,\qquad
Q_c=D_c\cap D_k,\qquad
Q_d=D_d\cap D_k.
\]

\begin{proposition}[Outer components]
\label{prop:outer-components-Ik-star}
Each of the outer components \(D_a,D_b,D_c,D_d\) contains at most one
specialization of an \(n\)-torsion section. Moreover, the torsion curve
\(T_n\) does not pass through any of the outer nodes
\(Q_a,Q_b,Q_c,Q_d\).
\end{proposition}

\begin{proof}
The outer components arise from resolving the quotient singularities at the
fixed points of \(\tau\) on \(C_0^{I_{2k}}\) and \(C_k^{I_{2k}}\). Each
outer component therefore corresponds to a single fixed point of \(\tau\)
on one of these two components.

Since these fixed points lie in the smooth locus of the \(I_{2k}\)-fiber, the \(n\)-torsion locus is étale over the base near them.
Consequently, at most one $n$-torsion section
can specialize to each fixed point. After passing to the quotient and
resolving, its strict transform meets the corresponding outer component at
a smooth point. Therefore each of \(D_a,D_b,D_c,D_d\) contains at most one
specialization of an \(n\)-torsion section.

Near a fixed point of $\tau$, choose analytic coordinates
$(z,s)$ in which
\(
\tau(z,s)=(-z,-s).
\)
The quotient is locally given by
\[
xt=y^2,\qquad x=z^2,\quad y=zs,\quad t=s^2.
\]
Resolving this $A_1$-singularity introduces the
corresponding outer component. The image of an invariant
torsion section meets this component away from its
intersection with the multiplicity-two component.

By the étaleness of the torsion locus in the smooth
locus of the semistable model, distinct torsion
sections have distinct specialization points.
Thus a torsion section specializing to a fixed
point of $\tau$ is invariant under $\tau$.
All remaining torsion specializations either lie at
smooth points exchanged by $\tau$, or at nodes
exchanged by $\tau$.
They therefore avoid the quotient singularities.
Consequently, $T_n$ does not pass through any of
the outer nodes.
\end{proof}

On the semistable \(I_{2k}\)-model,
Theorem~\ref{thm:torsion-distribution-Ik} gives
\[ \operatorname{length}(T_n^{I_{2k}}\cap (C_j^{I_{2k}})^*)
=
\begin{cases}
n,&2k\mid jn,\\
0,&2k\nmid jn.
\end{cases}
\]
To descend these points, we distinguish the fixed points of
\(\tau\) from the pairs that it exchanges.

\begin{theorem}[Smooth specialization for an \(I_k^*\)-fiber]
\label{thm:smooth-specialization-Ik-star}
The intersection lengths of \(T_n\) with the subsets \(D^*\) of the
reduced components are as follows.

\begin{enumerate}
\item At the left end,
\[
\operatorname{length}(T_n\cap D_a^*)=1,
\qquad
\operatorname{length}(T_n\cap D_b^*)=
\begin{cases}
1,&2\mid n,\\
0,&2\nmid n,
\end{cases}
\]
and
\[
\operatorname{length}(T_n\cap D_0^*)=
\begin{cases}
\dfrac{n-1}{2},&2\nmid n,\\[0.6em]
\dfrac{n-2}{2},&2\mid n.
\end{cases}
\]

\item At the right end,
\[
\operatorname{length}(T_n\cap D_c^*)
=
\operatorname{length}(T_n\cap D_d^*)
=
\begin{cases}
1,&2\mid n\text{ and }4\mid nk,\\
0,&\text{otherwise},
\end{cases}
\]
and
\[
\operatorname{length}(T_n\cap D_k^*)=
\begin{cases}
\dfrac{n-2}{2},
&2\mid n\text{ and }4\mid nk,\\[0.6em]
\dfrac n2,
&2\mid n\text{ and }4\nmid nk,\\[0.6em]
0,
&2\nmid n.
\end{cases}
\]

\item For \(1\leq j\leq k-1\),
\[
\operatorname{length}(T_n\cap D_j^*)=
\begin{cases}
n,&2k\mid jn,\\
0,&2k\nmid jn.
\end{cases}
\]
\end{enumerate}
\end{theorem}

\begin{proof}
The component \(C_0^{I_{2k}}\) contains \(n\) limiting \(n\)-torsion
points. The involution on these points is \(x\mapsto-x\). The zero point is fixed by $\tau$ and gives the zero section, which meets $D_a^*$.
If $n$ is even, the other fixed torsion point
gives a section meeting $D_b^*$. By Proposition~\ref{prop:outer-components-Ik-star},
these are the only possible torsion specializations on the two outer
components.

The remaining torsion points on
$C_0^{I_{2k}}$ occur in pairs exchanged
by $\tau$. Each pair gives one irreducible
component of $T_n$ of degree two over
$\Delta_t$, meeting $D_0^*$ once.
Hence
\[
\operatorname{length}(T_n\cap D_0^*)=
\begin{cases}
(n-1)/2,&2\nmid n,\\
(n-2)/2,&2\mid n,
\end{cases}
\]
which proves~(1).

The component \(C_k^{I_{2k}}\) contains limiting \(n\)-torsion if and only
if \(2k\mid kn\), equivalently if \(n\) is even. In that case it contains
\(n\) torsion points. The involution has two fixed torsion points on
\(C_k^{I_{2k}}\) precisely when \(4\mid nk\). These descend to sections
meeting \(D_c^*\) and \(D_d^*\), respectively. If \(4\nmid nk\), there are
no fixed torsion points on \(C_k^{I_{2k}}\).

The remaining torsion points occur in \(\tau\)-pairs and descend to
degree-two components meeting \(D_k\) once. Thus there are \((n-2)/2\)
such intersections when \(4\mid nk\), and \(n/2\) when \(4\nmid nk\).
This proves~(2).

Finally, let $1\leq j\leq k-1$. The map $\tau$
exchanges $C_j^{I_{2k}}$ and $C_{2k-j}^{I_{2k}}$,
which have the same image $D_j$ in the quotient.
These components contain limiting $n$-torsion
precisely when $2k\mid jn$. If this condition
holds, each contains exactly $n$ torsion points.
The map $\tau$ pairs the corresponding torsion
sections, one on each component. Each pair gives
one irreducible component of $T_n$ of degree two
over $\Delta$, meeting $D_j^*$ once. Therefore
\[
\operatorname{length}(T_n\cap D_j^*)=n.
\]
If \(2k\nmid jn\), neither component contains limiting \(n\)-torsion, so
\(\operatorname{length}(T_n\cap D_j^*)=0\).
\end{proof}

\subsection{Irreducible Components at the Central-Chain Nodes}
\label{subsec:central-nodes-Ik-star}

Under the quadratic semistable reduction, \(P_j\) has two preimages on the
\(I_{2k}\)-fiber, namely the nodes between
\(C_{j-1}^{I_{2k}}\) and \(C_j^{I_{2k}}\), and between
\(C_{2k-j}^{I_{2k}}\) and \(C_{2k-j+1}^{I_{2k}}\). These two nodes are exchanged by \(\tau\).

\begin{theorem}[Irreducible components at the central-chain nodes]
\label{thm:branches-Ik-star}
Fix \(1\leq j\leq k\). Then:

\begin{enumerate}
\item The node \(P_j\) is met by \(T_n\) precisely for those
\(i\in\{1,\ldots,n-1\}\) satisfying
\[
(j-1)n<2ki<jn.
\]

\item For every such $i$, set
$d_i:=\gcd(2ki,n)$.
There are $d_i$ irreducible components
\[
    B_{i,1},\ldots,B_{i,d_i}
\]
of $T_n$ specializing to $P_j$ and ramification index \(2n/d_i\) over \(0\).

\item The scheme-theoretic length at \(P_j\) is
\[
\operatorname{length}_{P_j}(T_n\cap\mathcal E_0)
=
2n\,
\#\left\{
1\leq i\leq n-1:
(j-1)n<2ki<jn
\right\}.
\]

\item For every \(i\) as above and \(1\leq c \leq d_i\),
\[
B_{i,c}\cdot D_{j-1}
=
\frac{jn-2ki}{d_i},
\qquad
B_{i, c}\cdot D_j
=
\frac{2ki-(j-1)n}{d_i},
\]
and \(B_{i,c}\) is disjoint from every other component of the central
fiber.
\end{enumerate}
\end{theorem}

\begin{proof}
Let $Q_j$ be the node between
$C_{j-1}^{I_{2k}}$ and $C_j^{I_{2k}}$.
The map $\tau$ exchanges $Q_j$ with the
other node lying above $P_j$.

By Theorem~\ref{thm:torsion-distribution-Ik},
applied to the $I_{2k}$-fiber, the irreducible
components of $T_n^{I_{2k}}$ specializing to
$Q_j$ are indexed by the integers
$i\in\{1,\ldots,n-1\}$ satisfying
\(
    (j-1)n<2ki<jn.
\)
For each such $i$, put $d_i=\gcd(2ki,n)$.
There are $d_i$ irreducible components
$B'_{i,1},\ldots,B'_{i,d_i}$ specializing
to $Q_j$, each of degree $n/d_i$ over $\Delta$.

The map $\tau$ sends these components to
the corresponding components specializing
to the other node above $P_j$. Each pair
has the same image in the original model
and gives one irreducible component
$B_{i,c}\subset T_n$ specializing to
$P_j$. Since $t=s^2$, the degree of
$B_{i,c}$ over $\Delta_t$ is $2n/d_i$.
By Proposition~\ref{prop:monodromy-components},
it has a unique analytic branch at $P_j$,
with ramification index $2n/d_i$.
This proves~(1) and~(2).

For each admissible $i$, the $d_i$
irreducible components contribute
\[
    d_i\cdot\frac{2n}{d_i}=2n
\]
to the intersection length at $P_j$.
Summing over all such $i$ gives
\[
    \operatorname{length}_{P_j}
    (T_n\cap\mathcal E_0)
    =
    2n\,\#\left\{
    1\leq i\leq n-1:
    (j-1)n<2ki<jn
    \right\},
\]
proving~(3).

Finally, Theorem~\ref{thm:branchwise-intersections-Ik}
gives
\[
    B'_{i,c}\cdot C_{j-1}^{I_{2k}}
    =\frac{jn-2ki}{d_i},
    \qquad
    B'_{i,c}\cdot C_j^{I_{2k}}
    =\frac{2ki-(j-1)n}{d_i}.
\]
Since $\tau$ exchanges the two nodes
lying above $P_j$, it has no fixed point
near either node. The quotient map is
therefore a local isomorphism near $Q_j$,
and the local intersection multiplicities
are unchanged. Consequently,
\[
    B_{i,c}\cdot D_{j-1}
    =\frac{jn-2ki}{d_i},
    \qquad
    B_{i,c}\cdot D_j
    =\frac{2ki-(j-1)n}{d_i}.
\]
Since no other component of the central
fiber passes through $P_j$, this proves~(4).
\end{proof}

\subsection{Component-wise Intersection Numbers}
\label{subsec:componentwise-intersections-Ik-star}

The smooth specialization points and the central-chain nodes account for
all intersections of \(T_n\) with the reduced central fiber. We sum their
contributions to obtain \(T_n\cdot D_j\).

\begin{corollary}[Component-wise intersection numbers]
\label{cor:componentwise-intersections-Ik-star}

For \(1\leq r\leq k-1\),
\[
T_n\cdot D_r
=
\sum_{i\in\mathbb Z}
\bigl(n-|nr-2ki|\bigr)_+.
\]

At the two ends of the multiplicity-two chain,
\[
T_n\cdot D_0
=
\frac12
\sum_{i\in\mathbb Z}
\bigl(n-|2ki|\bigr)_+
-
\frac{\gcd(n,2)}{2},
\]
and
\[
T_n\cdot D_k
=
\frac12
\sum_{i\in\mathbb Z}
\bigl(n-|nk-2ki|\bigr)_+
-
\mathds{1}_{\,2\mid n}\mathds{1}_{\,4\mid nk}.
\]

For the four outer components,
\[
T_n\cdot D_\alpha
=
\operatorname{length}(T_n\cap D_\alpha^*),
\qquad
\alpha\in\{a,b,c,d\}.
\]
\end{corollary}

\begin{proof}
We first consider an interior component \(D_r\), with
\(1\leq r\leq k-1\). By
Theorem~\ref{thm:branches-Ik-star}, a nodal contribution indexed by \(i\)
meets \(D_r\) precisely when
\(0<|nr-2ki|<n\), and the sum of the contributions of the corresponding
\(d_i\) branches is
\(
n-|nr-2ki|.
\)
If \(nr=2ki\), the torsion specializes to the smooth locus of
\(D_r\). By Theorem~\ref{thm:smooth-specialization-Ik-star}, this occurs
precisely when \(2k\mid rn\), and its contribution is \(n\). Thus the
smooth and nodal contributions are combined by the single formula
\[
T_n\cdot D_r
=
\sum_{i\in\mathbb Z}
\bigl(n-|nr-2ki|\bigr)_+.
\]

We next consider \(D_0\). Its nodal contributions are indexed by positive
integers \(i\) with \(2ki<n\), and their total contribution is
\[
\sum_{\substack{i\geq1\\2ki<n}}(n-2ki).
\]
Since the summand \((n-|2ki|)_+\) is symmetric under \(i\mapsto-i\), we have
\[
\frac12
\sum_{i\in\mathbb Z}
\bigl(n-|2ki|\bigr)_+
=
\frac n2
+
\sum_{\substack{i\geq1\\2ki<n}}(n-2ki).
\]
Here the term \(n/2\) comes from the central index \(i=0\). This index
corresponds to torsion specializing to the smooth locus at the left end.
By Theorem~\ref{thm:smooth-specialization-Ik-star}, the contribution on
\(D_0^*\) is
\[
\operatorname{length}(T_n\cap D_0^*)
=
\frac{n-\gcd(n,2)}{2}.
\]
Therefore the difference between the central term \(n/2\) in the symmetric
sum and the actual contribution on \(D_0^*\) is \(\gcd(n,2)/2\), and hence
\[
T_n\cdot D_0
=
\frac12
\sum_{i\in\mathbb Z}
\bigl(n-|2ki|\bigr)_+
-
\frac{\gcd(n,2)}{2}.
\]

Finally, consider \(D_k\). The summand
\(
\bigl(n-|nk-2ki|\bigr)_+
\)
is symmetric under \(i\mapsto n-i\). Hence
\[
\frac12
\sum_{i\in\mathbb Z}
\bigl(n-|nk-2ki|\bigr)_+
\]
accounts for one representative from each symmetric pair of nodal
contributions.

If \(n\) is even, there is also a central index \(i=n/2\), for which
\(nk-2ki=0\). Its contribution to the half-sum is \(n/2\), and this index
corresponds to torsion specializing to the smooth locus at the right end.

By Theorem~\ref{thm:smooth-specialization-Ik-star},
\[
\operatorname{length}(T_n\cap D_k^*)
=
\frac n2\,\mathds{1}_{\,2\mid n}
-
\mathds{1}_{\,2\mid n}\mathds{1}_{\,4\mid nk}.
\]

Thus the contribution of the central index to the half-sum differs from
the actual smooth-locus contribution by
\[
\mathds{1}_{\,2\mid n}\mathds{1}_{\,4\mid nk}.
\]
It follows that
\[
T_n\cdot D_k
=
\frac12
\sum_{i\in\mathbb Z}
\bigl(n-|nk-2ki|\bigr)_+
-
\mathds{1}_{\,2\mid n}\mathds{1}_{\,4\mid nk}.
\]
\end{proof}

\subsection{Connected Components}
\label{subsec:connected-components-Ik-star}

Every irreducible component of \(T_n\) has a unique specialization point on
\(\mathcal E_0\). Since distinct irreducible components are disjoint over
\(\Delta^*\), two of them lie in the same connected component precisely when
they specialize to the same point. Hence
\[
\#\pi_0(T_n)
=
\#\operatorname{Supp}(T_n\cap\mathcal E_0).
\]

\begin{proposition}[Connected components of \(T_n\)]
\label{prop:connected-components-Ik-star}
Set
\(
g=\gcd(n,2k).
\)
Then
\[
\#\pi_0(T_n)
=
\frac{ng+\gcd(n,2)}{2}
+
\mathds{1}_{\,2\mid n}\mathds{1}_{\,4\mid nk}
+
\min\left\{k,\frac{n-g}{2}\right\}.
\]
\end{proposition}

\begin{proof}
We count separately the smooth specialization points and the central-chain
nodes.

On the left end, the \(n\) torsion points on \(C_0^{I_{2k}}\) are
identified under the involution \(x\mapsto-x\). Since this involution has
\(\gcd(n,2)\) fixed points, they give
\[
\frac{n+\gcd(n,2)}{2}
\]
distinct smooth specialization points.

If \(n\) is even, the component \(C_k^{I_{2k}}\) also contains \(n\)
torsion points. Their quotient gives
\[
\frac n2+\mathds{1}_{\,4\mid nk}
\]
smooth specialization points. If \(n\) is odd, there is no torsion
specializing to \(C_k^{I_{2k}}\).

For \(1\leq j\leq k-1\), the component \(D_j^*\) contains smooth torsion
specializations precisely when \(2k\mid jn\). Since
\(g=\gcd(n,2k)\), this is equivalent to
\[
\frac{2k}{g}\mid j.
\]
There are
\[
\left\lfloor\frac{g-1}{2}\right\rfloor
\]
such indices \(j\), and each contributes \(n\) distinct specialization
points.

If \(n\) is odd, then \(g\) is odd, so the total number of smooth
specialization points is
\[
\frac{n+1}{2}
+
n\frac{g-1}{2}
=
\frac{ng+1}{2}.
\]
If \(n\) is even, then \(g\) is even, and the total is
\[
\frac{n+2}{2}
+
\left(\frac n2+\mathds{1}_{\,4\mid nk}\right)
+
n\left(\frac g2-1\right)
=
\frac{ng}{2}
+
1
+
\mathds{1}_{\,4\mid nk}.
\]
Both cases are summarized by
\[
\frac{ng+\gcd(n,2)}{2}
+
\mathds{1}_{\,2\mid n}\mathds{1}_{\,4\mid nk}.
\]

It remains to count the central-chain nodes. By
Theorem~\ref{thm:branches-Ik-star}, the node \(P_j\) is met precisely when
the interval
\(
((j-1)n,jn)
\)
contains a multiple of \(2k\).

Suppose first that \(n>2k\). Since each interval has length greater than
\(2k\), it contains a multiple of \(2k\). Hence all \(k\) central-chain
nodes are met.

Now suppose that \(n\leq 2k\). Each interval contains at most one multiple
of \(2k\). The relevant multiples are \(2ki\) with
\(
1\leq i<\frac n2.
\)
Such a multiple lies on an endpoint of one of the intervals precisely when
\(
n\mid 2ki.
\)
Writing \(g=\gcd(n,2k)\), this is equivalent to
\(
\frac{n}{g}\mid i.
\)
Among the integers \(1\leq i<n/2\), the number that do not satisfy this
divisibility condition is
\(
\frac{n-g}{2}.
\)
Therefore the number of central-chain nodes met by \(T_n\) is
\[
\min\left\{k,\frac{n-g}{2}\right\}.
\]

Adding the smooth and nodal specialization points proves the formula.
\end{proof}

\section{Good Reduction after Ramified Base Change}
\label{sec:potential-good-reduction}

The potentially good Kodaira fibers are
\[
II,\ III,\ IV,\ I_0^*,\ IV^*,\ III^*,\ II^*.
\]
Their local monodromy is finite, and after a suitable
ramified base change, followed by normalization,
resolution, and relative minimalization, the family
acquires good reduction
\cite[Lecture VI, \S 4, Table VI.4.1]{Miranda1989}.
The required base-change degrees are
\[
\begin{array}{c|ccccccc}
\text{fiber type}
& II & III & IV & I_0^* & IV^* & III^* & II^* \\ \hline
\text{degree}
& 6 & 4 & 3 & 2 & 3 & 4 & 6.
\end{array}
\]

Throughout this section, by the relatively minimal model of a ramified
base change we mean the normalization of the fiber product, its minimal
resolution, and the contraction of vertical $(-1)$-curves. In several cases,
a component of the central fiber gives rise to a cyclic cover of
$\mathbf P^1$, whose ramification is determined by the multiplicities of
adjacent components. We use Riemann--Hurwitz to identify its genus.

In the next section, we study the action induced by
$s\mapsto\zeta_h s$ on the good-reduction model after the base change
$t=s^h$ and use it to determine the specialization of $n$-torsion points.

\subsection{The Fibers of Type
\texorpdfstring{\(II\), \(III\), and \(IV\)}{II, III, and IV}}
\label{subsec:potential-good-III-IV}

A fiber of type $II$ is an irreducible cuspidal rational curve. Resolving
the cusp produces a component of multiplicity $6$, which gives the smooth
central fiber after a degree-six base change.

\begin{proposition}[Good reduction for a type $II$ fiber]
\label{prop:II}
Let $\pi\colon\mathcal E\to\Delta$ have central fiber of type $II$.
After a totally ramified base change of degree $6$, the relatively minimal
model has smooth central fiber.
\end{proposition}

\begin{proof}
Resolving the cusp by three successive blow-ups produces components of
multiplicities $2$, $3$, and $6$. The last exceptional curve
$E_3\simeq\mathbf P^1$ meets components of multiplicities $1$, $2$,
and $3$ at three distinct points.

After the degree-six base change, the normalized inverse image of $E_3$
is a connected cyclic cover
\[
    \widetilde E\longrightarrow E_3\simeq\mathbf P^1
\]
of degree $6$. The ramification indices over these three points are
$6$, $3$, and $2$, respectively. The Riemann--Hurwitz formula gives
\[
    2g(\widetilde E)-2
    =6(-2)+(6-1)+2(3-1)+3(2-1)=0.
\]
Thus $\widetilde E$ has genus one. The remaining components are six
rational tails of multiplicity one, each meeting $\widetilde E$
transversely once. Since the intersection of each tail with the central
fiber is zero, each has self-intersection $-1$. Contracting the tails
leaves $\widetilde E$ as the smooth central fiber.
\end{proof}

For types $III$ and $IV$, the central fibers have a non-normal-crossing
point. Resolving it introduces a component of multiplicity $4$ in type
$III$ and a component of multiplicity $3$ in type $IV$. We use these
components in the corresponding base changes; see
\cite{Kodaira1963,Miranda1989}.

\begin{proposition}[Potentially good reduction for a type \(III\) fiber]
\label{prop:III}
Let $\pi\colon\mathcal E\to\Delta$ be a relatively minimal elliptic fibration whose
central fiber has type \(III\). After two blow-ups at the
non-normal-crossing point and a totally ramified base change of degree \(4\),
the relatively minimal pullback has smooth central fiber.
\end{proposition}

\begin{proof}
Write
\(
F=C_1+C_2,
\)
where \(C_1\) and \(C_2\) are smooth rational curves tangent at one point.
Blowing up the tangency point gives
\[
F_1=C_1'+C_2'+2E_1.
\]
The three components meet at one point. Blowing up this point gives
\[
F_2=C_1''+C_2''+2E_1'+4E_2.
\]
The exceptional curve $E_2\simeq\mathbf P^1$ has multiplicity $4$ and
meets \(C_1''\), \(C_2''\), and \(E_1'\) transversely at three distinct
points.

After the degree-four base change, the normalized inverse image of \(E_2\)
is a cyclic cover
\[
\widetilde E\longrightarrow E_2\simeq\mathbf P^1.
\]
The adjacent components have multiplicities $1$, $1$, and $2$, so the
ramification indices are $4$, $4$, and $2$, respectively. Above the last
point there are two ramification points of index $2$. Thus
Riemann--Hurwitz gives
\[
2g(\widetilde E)-2
=
4(-2)+(4-1)+(4-1)+2
=
0.
\]
Hence \(\widetilde E\) is a smooth elliptic curve.

On the regular pullback before relative minimalization, the remaining
components form four rational tails
\[
F'=\widetilde E+\sum_{i=1}^4T_i,
\]
each meeting \(\widetilde E\) transversely once. Since \(F'\cdot T_i=0\),
\(
T_i^2=-1.
\)
Contracting these four tails leaves \(\widetilde E\) as the smooth central
fiber.
\end{proof}

For type \(IV\), one blow-up already produces the component of multiplicity three.

\begin{proposition}[Potentially good reduction for a type \(IV\) fiber]
\label{prop:IV}
Let $\pi\colon\mathcal E\to\Delta$ be a relatively minimal elliptic fibration whose
central fiber has type \(IV\). After blowing up the triple point and
performing a totally ramified base change of degree \(3\), the relatively
minimal pullback has smooth central fiber.
\end{proposition}

\begin{proof}
Write
\(
F=C_1+C_2+C_3,
\)
where the three smooth rational components meet at one point. Blowing up this
point gives
\[
F_1=C_1'+C_2'+C_3'+3E,
\]
where \(E\simeq\mathbf P^1\) and the strict transforms \(C_i'\) are disjoint, each meeting \(E\) transversely once.

After the cubic base change, the normalized inverse image of \(E\) is a
cyclic triple cover
\[
\widetilde E\longrightarrow E\simeq\mathbf P^1
\]
totally ramified at the three points \(E\cap C_i'\). Hence
\[
2g(\widetilde E)-2
=
3(-2)+3(3-1)
=
0,
\]
so \(\widetilde E\) is a smooth elliptic curve.

Before relative minimalization, the remaining components form three rational tails
\[
F'=\widetilde E+T_1+T_2+T_3,
\]
each meeting \(\widetilde E\) transversely once. Since \(F'\cdot T_i=0\),
we have \(T_i^2=-1\). Contracting the three tails leaves
\(\widetilde E\) as the smooth central fiber.
\end{proof}

\subsection{The Starred Fibers}
\label{subsec:potential-good-starred}

For the starred fibers, the component of maximal multiplicity is already
present in the Kodaira configuration. We treat $I_0^*$ and $IV^*$ first;
their good-reduction models will also be used for $III^*$ and $II^*$.

For a fiber of type $I_0^*$, write
\[
    \mathcal E_0=C_0+C_1+C_2+C_3+2C_4,
\]
where $C_4\simeq\mathbf P^1$ is the central
multiplicity-two component and the four
multiplicity-one components $C_0,\ldots,C_3$
meet $C_4$ transversely at distinct points.
We label them so that the zero section meets
$C_0$. These four intersection points form
the branch locus of the double cover obtained
after the quadratic base change.

\begin{proposition}[Potentially good reduction for a type $I_0^*$ fiber]
\label{prop:I0star}
Let $\pi\colon\mathcal E\to\Delta$ be a relatively
minimal elliptic fibration whose central fiber
has type $I_0^*$. A totally ramified quadratic
base change produces good reduction after
normalization, resolution, and relative
minimalization.
\end{proposition}

\begin{proof}
After the quadratic base change, the normalized
inverse image of $C_4$ is the double cover
\[
    E\longrightarrow C_4\simeq\mathbf P^1
\]
branched at the four points $C_4\cap C_i$,
$0\leq i\leq3$. Riemann--Hurwitz gives
\[
    2g(E)-2=2(-2)+4=0,
\]
so $E$ is a smooth elliptic curve.

Before relative minimalization, the remaining
components form four rational curves $R_i$,
$0\leq i\leq3$, each meeting $E$ transversely
once. Writing the central fiber of this model
as $F'=E+\sum_{i=0}^3R_i$, we have
\[
    0=F'\cdot R_i=R_i^2+1.
\]
Thus each $R_i$ is a $(-1)$-curve.
Contracting these four curves leaves $E$
as the smooth central fiber.
\end{proof}

For type $IV^*$, write
$F=3C_0+\sum_{i=1}^3(2C_i+D_i)$, where each $C_i$ meets $C_0$
and $D_i$ meets $C_i$. The multiplicity-three component $C_0$
gives a triple cover branched at the three points $C_0\cap C_i$.

\begin{proposition}[Potentially good reduction for a type \(IV^*\) fiber]
\label{prop:IVstar}
Let $\pi\colon\mathcal E\to\Delta$ be a relatively minimal elliptic fibration whose
central fiber has type \(IV^*\). A totally ramified base change of degree \(3\)
produces good reduction after normalization, resolution, and relative
minimalization.
\end{proposition}

\begin{proof}
After the cubic base change, the normalized inverse image of \(C_0\) is a
cyclic triple cover
\[
E\longrightarrow C_0\simeq\mathbf P^1
\]
totally ramified at the three points \(C_0\cap C_i\). Hence
Riemann--Hurwitz gives
\[
2g(E)-2
=
3(-2)+3(3-1)
=
0,
\]
so \(E\) is a smooth elliptic curve.

After normalization, resolution, and contraction of the
outermost $(-1)$-curves, we obtain a regular model whose
special fiber has the form
\[
F'
=
E+\sum_{i=1}^3(2R_i+T_i),
\]
where each \(R_i\) meets \(E\) and \(T_i\) transversely. Since
\(F'\cdot R_i=0\),
\[
0=2R_i^2+1+1,
\]
so \(R_i^2=-1\). Similarly,
\[
0=F'\cdot T_i=T_i^2+2,
\]
and hence \(T_i^2=-2\).

After contracting \(R_i\), the image of \(T_i\) has self-intersection
\(-1\), so it can also be contracted. Repeating this for the three arms
leaves \(E\) as the smooth central fiber.
\end{proof}

The fibers of types $III^*$ and $II^*$ can be treated by successive
base changes through the preceding models:
\[
III^*
\xrightarrow{\;2\;}
I_0^*
\xrightarrow{\;2\;}
I_0,
\qquad
II^*
\xrightarrow{\;2\;}
IV^*
\xrightarrow{\;3\;}
I_0.
\]
Each arrow denotes a totally ramified base change
of the indicated degree, followed by normalization,
resolution, and relative minimalization
\cite[Lecture VI, \S 4, Table VI.4.1]{Miranda1989}.

\begin{proposition}[Potentially good reduction for a type \(III^*\) fiber]
\label{prop:IIIstar}
Let $\pi\colon\mathcal E\to\Delta$ be a relatively minimal elliptic fibration whose
central fiber has type \(III^*\). Two successive totally ramified quadratic
base changes produce good reduction after normalization, resolution, and
relative minimalization. Equivalently, the total ramification degree is \(4\).
\end{proposition}

\begin{proof}
Under the first quadratic base change \(t=s^2\), the local monodromy is
squared. For a fiber of type \(III^*\), the monodromy has order \(4\), so
after this base change it becomes \(-I\), which is the local monodromy of a
fiber of type \(I_0^*\). After normalization, resolution, and relative
minimalization, the resulting central fiber therefore has type \(I_0^*\). 

Applying Proposition~\ref{prop:I0star} to this intermediate \(I_0^*\)-fiber,
a second quadratic base change produces good reduction. Thus the composite
base change has total ramification degree \(4\).
\end{proof}

\begin{proposition}[Potentially good reduction for a type \(II^*\) fiber]
\label{prop:IIstar}
Let $\pi\colon\mathcal E\to\Delta$ be a relatively minimal elliptic fibration whose
central fiber has type \(II^*\). A totally ramified quadratic base change
followed by a totally ramified cubic base change produces good reduction after
normalization, resolution, and relative minimalization. Equivalently, the
total ramification degree is \(6\).
\end{proposition}

\begin{proof}
Under the first quadratic base change \(t=s^2\), the local monodromy is
squared. For a fiber of type \(II^*\), the monodromy has order \(6\), so
after this base change it has order $3$. The Kodaira base-change table
identifies the resulting relatively minimal fiber as type $IV^*$.

Applying the cubic good-reduction construction of
Proposition~\ref{prop:IVstar} to this intermediate \(IV^*\)-fiber gives a
smooth central elliptic curve. Thus the composite base change has total
ramification degree \(6\) and produces good reduction.
\end{proof}

\section{Torsion Curves over Fibers with Potentially Good Reduction}
\label{sec:additive-torsion}

We use the good-reduction models of
Section~\ref{sec:potential-good-reduction} to study the specialization of
$T_n$ for fibers of type $II$, $III$, $IV$, $I_0^*$, $IV^*$, $III^*$,
and $II^*$. After a base change $t=s^h$ giving good reduction, the
$n$-torsion points extend to sections of the smooth elliptic model.
Here $\zeta_h$ denotes a primitive $h$th root of unity. The change
$s\mapsto\zeta_hs$ acts on these sections. Comparing their
specialization points and their images on the original surface gives the
specialization of the irreducible components of $T_n$.

For any irreducible component $H$ of the central fiber, we write
$H^*:=H\setminus\operatorname{Sing}((\mathcal E_0)_{\mathrm{red}})$
and denote by $m_H$ its multiplicity in $\mathcal E_0$.
If an irreducible component $B\subset T_n$ of degree $r$ over $\Delta$
specializes to a point of $H^*$, then
\[
    r=B\cdot\mathcal E_0=m_H(B\cdot H),
    \qquad B\cdot H=\frac{r}{m_H}.
\]
In particular, if $m_H=r$, then $B$ meets $H$ transversely.
The distinction between multiplicity one and higher multiplicities
will be useful for the starred fibers.

\subsection{Finite Monodromy and Orbit Lengths}
\label{subsec:finite-additive-monodromy}

Identify $E_t[n]\simeq(\mathbb Z/n\mathbb Z)^2$,
and let $M$ denote the integral local monodromy
matrix, acting on $E_t[n]$ through its reduction
modulo $n$. $M$ has order $6$, $4$, $3$, or $2$ for fibers of
type $II$, $III$, $IV$, or $I_0^*$, respectively. For the paired types
$(II,II^*)$, $(III,III^*)$, and $(IV,IV^*)$, the monodromy matrices
are inverse up to conjugacy and hence have the same orbits on $E_t[n]$. 

\begin{proposition}[Orbit structure for finite additive monodromy]
\label{prop:finite-monodromy-orbits}
Identify $E_t[n]\simeq(\mathbb Z/n\mathbb Z)^2$ and let $M$
denote the local monodromy matrix modulo $n$.

\begin{enumerate}
    \item For fibers of type \(II\) or \(II^*\), the integral monodromy has order \(6\).
    The origin is the only fixed point. If \(3\mid n\), there is one orbit
    of length \(2\); if \(2\mid n\), there is one orbit of length \(3\).
    Every remaining nonzero orbit has length \(6\).

    \item For fibers of type \(III\) or \(III^*\), the integral monodromy has order \(4\).
    If \(n\) is odd, the origin is the only fixed point and every nonzero
    orbit has length \(4\). If \(n\) is even, there are two fixed points,
    one orbit of length \(2\), and every remaining orbit has length \(4\).

    \item For fibers of type \(IV\) or \(IV^*\), the integral monodromy has order \(3\).
    It has one fixed point if \(3\nmid n\) and three fixed points if
    \(3\mid n\). Every non-fixed orbit has length \(3\).

    \item For a fiber of type \(I_0^*\), \(M=-I\). If \(n\) is odd, the
    origin is the only fixed point and every nonzero orbit has length \(2\).
    If \(n\) is even, the four points of \(E_t[2]\) are fixed and
    every remaining orbit has length \(2\).
\end{enumerate}

In every case, the orbit length is the degree over \(\Delta\) of the
corresponding irreducible component of \(T_n\).
\end{proposition}

\begin{proof}
After choosing suitable bases, we may take the finite-order monodromy matrices
to be
\[
M_6=
\begin{pmatrix}
1&1\\
-1&0
\end{pmatrix},
\qquad
M_4=
\begin{pmatrix}
0&1\\
-1&0
\end{pmatrix},
\qquad
M_3=
\begin{pmatrix}
0&1\\
-1&-1
\end{pmatrix},
\qquad
M_2=-I,
\]
of orders \(6,4,3,2\), respectively; see
\cite{Kodaira1963,Miranda1989}. The starred duals have inverse monodromy,
which gives the same orbits.

For \(M_6\), the matrix \(M_6-I\) is invertible over
\(\mathbb Z/n\mathbb Z\), so the origin is the only fixed point. Moreover,
\(
M_6^3=-I,
\)
while \(M_6^2\) has order \(3\). If \(3\mid n\), the two nonzero points fixed
by \(M_6^2\) form one orbit of length \(2\). If \(2\mid n\), the three
nonzero points fixed by \(M_6^3=-I\) form one orbit of length \(3\).
All remaining points have trivial stabilizer and hence orbit length \(6\).

For \(M_4\), one has
\(
M_4^2=-I.
\)
If \(n\) is odd, \(-I\) fixes only the origin, so every nonzero orbit has
length \(4\). If \(n\) is even, there are four points satisfying \(2v=0\).
Two of them are fixed by \(M_4\), while the other two form a single orbit of
length \(2\). Every other point has orbit length \(4\).

For \(M_3\), solving
\(
(M_3-I)v=0
\)
gives one solution when \(3\nmid n\) and three solutions when \(3\mid n\).
Since \(M_3^3=I\), every non-fixed orbit has length \(3\).

Finally, for \(M_2=-I\), the fixed points are precisely the solutions of
\(
2v=0.
\)
There is only the origin when \(n\) is odd, while for even \(n\) these are
the four points of \(E_t[2]\). Every other point is paired with its negative and hence belongs to an orbit of length \(2\).

The final assertion follows from
Proposition~\ref{prop:monodromy-components}.
\end{proof}

\subsection{The Unstarred Fibers}
\label{subsec:unstarred-additive}

For types $II$, $III$, and $IV$, the points fixed by monodromy give
sections meeting the smooth locus of the central fiber. Every other
orbit specializes to its unique singular point. We obtain the
intersection length there from the orbit counts; for types $III$ and
$IV$, we also compute intersections with the individual fiber components.

\subsubsection*{Type \(II\)}

A type $II$ fiber is irreducible, with a single cusp.

\begin{theorem}[Torsion specialization for a type \(II\) fiber]
\label{thm:type-II}
Let the central fiber $C=\mathcal E_0$ have type $II$, and let $R$
be its cusp. Put $C^*=C\setminus\{R\}$.
Then:
\begin{enumerate}
    \item
    \[
    \operatorname{length}(T_n\cap C^*)=1
    \]

    \item
    \[
    \operatorname{length}_{R}(T_n\cap\mathcal E_0)=n^2-1.
    \]
    If \(2\mid n\), there is one irreducible component of degree \(3\) through \(R\);
    if \(3\mid n\), there is one irreducible component of degree \(2\) through \(R\);
    every remaining irreducible component specializing to $R$ has degree $6$.

    \item The number of irreducible components is
    \[
    \#\operatorname{Irr}(T_n)=
    \begin{cases}
    \dfrac{n^2-1}{6}+1,
        &2\nmid n\text{ and }3\nmid n,\\[0.8em]
    \dfrac{n^2-3}{6}+2,
        &3\mid n\text{ and }2\nmid n,\\[0.8em]
    \dfrac{n^2-4}{6}+2,
        &2\mid n\text{ and }3\nmid n,\\[0.8em]
    \dfrac{n^2-6}{6}+3,
        &6\mid n.
    \end{cases}
    \]
\end{enumerate}
\end{theorem}

\begin{proof}
The origin is the only fixed point of the monodromy, by
Proposition~\ref{prop:finite-monodromy-orbits}. It gives the zero section, which meets \(C^*\), while every
non-fixed orbit specializes to the cusp \(R\). Since \(T_n\to\Delta\) has
degree \(n^2\), this gives
\[
\operatorname{length}_{R}(T_n\cap\mathcal E_0)=n^2-1.
\]

The possible ramification degrees are the nontrivial orbit lengths
\(2,3,\) and \(6\). When \(3\mid n\), there is one orbit of length \(2\);
when \(2\mid n\), there is one orbit of length \(3\); all remaining
nonzero points lie in orbits of length \(6\). Counting these orbits together
with the zero section gives the stated number of irreducible components.
\end{proof}

\subsubsection*{Type \(III\)}

For a fiber of type $III$, write
$\mathcal E_0=C_1+C_2$, where $C_1$ meets
the zero section. The two components are
tangent at their unique intersection point
$R=C_1\cap C_2$, with intersection
multiplicity $2$. Put
$C_j^*:=C_j\setminus\{R\}$ for $j=1,2$.

\begin{theorem}[Torsion specialization for a type $III$ fiber]
\label{thm:type-III}
With the notation above, the following hold:
\begin{enumerate}
    \item
    \[
    \operatorname{length}(T_n\cap C_1^*) =1,
    \qquad
    \operatorname{length}(T_n\cap C_2^*) =
    \begin{cases}
    1,&2\mid n,\\
    0,&2\nmid n.
    \end{cases}
    \]

    \item
    \[
    \operatorname{length}_{R}(T_n\cap\mathcal E_0)=
    \begin{cases}
    n^2-2,&2\mid n,\\
    n^2-1,&2\nmid n.
    \end{cases}
    \]
    When \(2\mid n\), exactly one irreducible component specializing to \(R\) has degree \(2\),
    and every other irreducible component specializing to \(R\) has degree \(4\). When \(n\) is
    odd, every irreducible component specializing to \(R\) has degree \(4\).

    \item If $B$ is an irreducible component of $T_n$ specializing to $R$, then
    \[
    C_1\cdot B=C_2\cdot B=
    \begin{cases}
    1,&\deg(B/\Delta)=2,\\
    2,&\deg(B/\Delta)=4.
    \end{cases}
    \]

    \item The number of irreducible components is
    \[
    \#\operatorname{Irr}(T_n)=
    \begin{cases}
    \dfrac{n^2-4}{4}+3,&2\mid n,\\[0.8em]
    \dfrac{n^2-1}{4}+1,&2\nmid n.
    \end{cases}
    \]
\end{enumerate}
\end{theorem}

\begin{proof}

Proposition~\ref{prop:finite-monodromy-orbits} gives one fixed point when
\(n\) is odd and two when \(n\) is even. In the degree-four good-reduction
model of Proposition~\ref{prop:III}, the two components \(C_1\) and \(C_2\)
arise from rational tails attached to the smooth elliptic fiber at two
distinct fixed points of the action induced by $s\mapsto\zeta_4s$. Thus each of \(C_1^*\) and
\(C_2^*\) can contain the specialization of at most one fixed torsion point.

The origin gives the zero section meeting \(C_1^*\). Hence, when \(n\) is
even, the second fixed point specializes to \(C_2^*\). Every non-fixed orbit
specializes to the common singular point \(R\). Subtracting the smooth
contributions from the total degree \(n^2\) gives the stated length at \(R\).
The same proposition gives one orbit of length \(2\) when \(n\) is even,
with all remaining non-fixed orbits of length \(4\).

It remains to compute the local intersections at \(R\). Analytically locally, we may
choose coordinates in which
\[
C_1=\{y=0\},
\qquad
C_2=\{y-x^2=0\},
\qquad
t=y(y-x^2)
\]
up to multiplication of \(t\) by a unit. Hence, for any branch \(B\)
through \(R\),
\[
\deg(B/\Delta)
=
C_1\cdot B+C_2\cdot B.
\]

For a degree-two component, both intersection multiplicities are positive, so
they must both equal \(1\). For a degree-four component, neither intersection
multiplicity can equal \(1\): if
\(\operatorname{ord}_s y=1\), then
\(\operatorname{ord}_s(y-x^2)=1\) as well, and conversely. Since the two
positive intersection multiplicities sum to \(4\), they must therefore both
equal \(2\).

Finally, when \(n\) is even, the two fixed points and the two points in the
length-two orbit account for four torsion points, leaving
\[
\frac{n^2-4}{4}
\]
orbits of length \(4\). When \(n\) is odd, only the origin is fixed, leaving
\[
\frac{n^2-1}{4}
\]
orbits of length \(4\). This gives the irreducible-component count.
\end{proof}

\subsubsection*{Type \(IV\)}

For a fiber of type $IV$, write
$\mathcal E_0=C_1+C_2+C_3$, where $C_1$
meets the zero section. The three components
meet at a single point $R$, with distinct
tangent directions. Put
$C_j^*:=C_j\setminus\{R\}$ for $j=1,2,3$.

\begin{theorem}[Torsion specialization for a type $IV$ fiber]
\label{thm:type-IV}
With the notation above, the following hold:
\begin{enumerate}
    \item
    \[
    \operatorname{length}(T_n\cap C_1^*)=1,
    \qquad
    \operatorname{length}(T_n\cap C_2^*)= \operatorname{length}(T_n\cap C_3^*)=
    \begin{cases}
    1,&3\mid n,\\
    0,&3\nmid n.
    \end{cases}
    \]

    \item
    \[
    \operatorname{length}_{R}(T_n\cap\mathcal E_0)=
    \begin{cases}
    n^2-3,&3\mid n,\\
    n^2-1,&3\nmid n.
    \end{cases}
    \]
    Every irreducible component specializing to \(R\) has degree \(3\) over \(\Delta\).

    \item If $B$ is an irreducible component of $T_n$ specializing to $R$, then
    \[
    C_1\cdot B=C_2\cdot B=C_3\cdot B=1.
    \]

    \item The number of irreducible components is
    \[
    \#\operatorname{Irr}(T_n)=
    \begin{cases}
    \dfrac{n^2-3}{3}+3,&3\mid n,\\[0.8em]
    \dfrac{n^2-1}{3}+1,&3\nmid n.
    \end{cases}
    \]
\end{enumerate}
\end{theorem}
\begin{proof}
By Proposition~\ref{prop:finite-monodromy-orbits}, there is one fixed point
if $3\nmid n$ and three if $3\mid n$. In the cubic good-reduction model
of Proposition~\ref{prop:IV}, the three components $C_1,C_2,C_3$ arise
from rational tails at the fixed points of the action induced by
$s\mapsto\zeta_3s$. Thus each smooth locus
\(C_j^*\) can contain the specialization of at most one fixed torsion point.

The origin gives the zero section meeting \(C_1^*\). Hence, when \(3\mid n\),
the other two fixed points specialize to \(C_2^*\) and \(C_3^*\). Every
non-fixed orbit has length \(3\) and specializes to the common singular point
\(R\). Subtracting the smooth contributions from the total degree \(n^2\)
gives
\[
\operatorname{length}_{R}(T_n\cap\mathcal E_0)=
\begin{cases}
n^2-3,&3\mid n,\\
n^2-1,&3\nmid n.
\end{cases}
\]

As in the proof of Theorem~\ref{thm:type-III}, the degree of an
irreducible component specializing to $R$ is the sum of its
intersection multiplicities with the reduced components through $R$.
Thus, for such a component $B$,
\[
3=\deg(B/\Delta)
=C_1\cdot B+C_2\cdot B+C_3\cdot B.
\]
All three intersection multiplicities are positive, so
\[
C_1\cdot B=C_2\cdot B=C_3\cdot B=1.
\]

Finally, the non-fixed torsion points form orbits of length \(3\). Therefore
there are
\[
\frac{n^2-3}{3}
\]
such orbits when \(3\mid n\), and
\[
\frac{n^2-1}{3}
\]
when \(3\nmid n\). Adding the fixed sections gives the stated number of
irreducible components.
\end{proof}

\subsection{The Starred Fibers}
\label{subsec:starred-additive}

For $I_0^*$, $IV^*$, $III^*$, and $II^*$, an orbit of length $r$
specializes to a point of $H^*$ on a fiber component $H$ of
multiplicity $r$. By the intersection formula at the beginning of
the section, the corresponding irreducible component of $T_n$
meets $H$ with multiplicity one. Different orbits give different
specialization points. These facts also imply that the irreducible
components of $T_n$ are pairwise disjoint for these fiber types.

The following table records the component met by each orbit length,
using the notation of Section~\ref{sec:potential-good-reduction} for
$I_0^*$ and $IV^*$ and the labels specified below for $III^*$ and $II^*$.
\[
\begin{array}{c|llll}
\text{fiber type}&\multicolumn{4}{c}{\text{orbit length }\longmapsto\text{ component}}\\ \hline
I_0^* & 1\mapsto C_0,\ldots,C_3
      & 2\mapsto C_4 & & \\
IV^*&1\mapsto D_1,D_2,D_3&3\mapsto C_0&&\\
III^*&1\mapsto C_0,C_1&2\mapsto C_2&4\mapsto C_3&\\
II^*&1\mapsto C_0&2\mapsto C_1&3\mapsto C_2&6\mapsto C_3
\end{array}
\]

\subsubsection*{Type $I_0^*$}

For a fiber of type $I_0^*$, write
\[
    \mathcal E_0=C_0+C_1+C_2+C_3+2C_4,
\]
where $C_4$ is the central multiplicity-two component
and $C_0,\ldots,C_3$ have multiplicity one.
We label the components so that the zero section
meets $C_0$.

By Proposition~\ref{prop:I0star}, the quadratic
base change gives a double cover $E\to C_4$,
branched at the four points $C_4\cap C_i$,
$0\leq i\leq3$. These correspond to the four
fixed points of the map $P\mapsto-P$ on $E$.
Put $C_i^*:=C_i\setminus
\operatorname{Sing}((E_0)_{\mathrm{red}})$.

\begin{theorem}[Torsion specialization for a type $I_0^*$ fiber]
\label{thm:type-I0star}
With the notation above, the following hold:
\begin{enumerate}
    \item If $n$ is odd,
    \[
        \operatorname{length}(T_n\cap C_0^*)=1,
        \qquad
        \operatorname{length}(T_n\cap C_1^*)
        = \operatorname{length}(T_n\cap C_2^*)
        =\operatorname{length}(T_n\cap C_3^*)=0.
    \]
    If $n$ is even,
    \[
        \operatorname{length}(T_n\cap C_j^*)=1,
        \qquad 0\leq j\leq3.
    \]

    \item No singular point of $\mathcal E_0$
    lies on $T_n$.

    \item
    \[
    \operatorname{length}(T_n\cap C_4^*)=
        \begin{cases}
            \dfrac{n^2-1}{2},&2\nmid n,\\[0.8em]
            \dfrac{n^2-4}{2},&2\mid n.
        \end{cases}
    \]
    Every irreducible component of $T_n$ meeting
    $C_4^*$ has degree $2$ over $\Delta$,
    intersects $C_4$ with multiplicity one,
    and distinct such components meet $C_4^*$
    at distinct points.

    \item
    \[
        \#\operatorname{Irr}(T_n)=
        \begin{cases}
            \dfrac{n^2-1}{2}+1,&2\nmid n,\\[0.8em]
            \dfrac{n^2-4}{2}+4,&2\mid n.
        \end{cases}
    \]
\end{enumerate}
\end{theorem}

\begin{proof}
By Proposition~\ref{prop:finite-monodromy-orbits},
the monodromy is $-I$. Thus the origin is the
only fixed torsion point when $n$ is odd, while
there are four fixed points when $n$ is even;
every other orbit has length $2$.

The four multiplicity-one components correspond
to distinct fixed points of $P\mapsto-P$ on $E$.
Hence each $C_i^*$, $0\leq i\leq3$, contains
at most one specialization of an $n$-torsion
section.

The origin gives the zero section meeting $C_0^*$.
If $n$ is even, all four fixed points are
$n$-torsion, so the remaining three specialize
to $C_1^*,C_2^*,C_3^*$. If $n$ is odd, no other
fixed torsion point occurs.

The remaining torsion points lie in orbits of
length $2$ and specialize to points of $C_4^*$.
The fixed sections meet the smooth loci of the
outer components. Thus $T_n$ contains no
singular point of the central fiber.

Every non-fixed orbit has length $2$.
If $B_O\subset T_n$ is the corresponding
irreducible component, then
\[
    2=B_O\cdot\mathcal E_0
     =2(B_O\cdot C_4),
    \qquad B_O\cdot C_4=1.
\]
Distinct orbits have distinct images in $C_4^*$.

There are $n^2-1$ non-fixed torsion points
when $n$ is odd and $n^2-4$ when $n$ is even.
Dividing by the orbit length $2$ gives the
stated value of $C_4^*\cdot T_n$, and adding
the fixed sections gives the number of
irreducible components.
\end{proof}

\subsubsection*{Type \(IV^*\)}

For a fiber of type $IV^*$, write
\[
    \mathcal E_0
    =3C_0+\sum_{i=1}^3(2C_i+D_i),
\]
where $C_0$ is the multiplicity-three component,
each $C_i$ meets $C_0$, and $D_i$ meets $C_i$.
We label the components so that the zero
section meets $D_1$.

By Proposition~\ref{prop:IVstar}, the cubic
base change gives a triple cover $E\to C_0$,
totally ramified over the three points
$C_0\cap C_i$. These points correspond to
the three arms of the original fiber,
ending in $D_1,D_2,D_3$.

\begin{theorem}[Torsion specialization for a type $IV^*$ fiber]
\label{thm:type-IVstar}
With the notation above, the following hold:
\begin{enumerate}
    \item
    \[
    \operatorname{length}(T_n\cap D_1^*)=1,
    \qquad
    \operatorname{length}(T_n\cap D_2^*)=\operatorname{length}(T_n\cap D_3^*)=
    \begin{cases}
    1,&3\mid n,\\
    0,&3\nmid n.
    \end{cases}
    \]

    \item No singular point of \(\mathcal E_0\) lies on \(T_n\).

    \item
    \[
    \operatorname{length}(T_n\cap C_1^*)=\operatorname{length}(T_n\cap C_2^*)=\operatorname{length}(T_n\cap C_3^*)=
    0.
    \]

    \item
    \[
   \operatorname{length}(T_n\cap C_0^*)=
    \begin{cases}
    \dfrac{n^2-3}{3},&3\mid n,\\[0.8em]
    \dfrac{n^2-1}{3},&3\nmid n.
    \end{cases}
    \]
    Every component meeting \(C_0^*\) has degree \(3\) over \(\Delta\),
    intersects the reduced component \(C_0\) with multiplicity one, and
    distinct such components meet \(C_0^*\) at distinct points.

    \item
    \[
    \#\operatorname{Irr}(T_n)=
    \begin{cases}
    \dfrac{n^2-3}{3}+3,&3\mid n,\\[0.8em]
    \dfrac{n^2-1}{3}+1,&3\nmid n.
    \end{cases}
    \]
\end{enumerate}
\end{theorem}

\begin{proof}
By Proposition~\ref{prop:finite-monodromy-orbits}, there is one fixed point
if \(3\nmid n\) and three fixed points if \(3\mid n\); every non-fixed orbit
has length \(3\).

In the cubic good-reduction model of Proposition~\ref{prop:IVstar}, the
multiplicity-three component \(C_0\) is the quotient of the smooth elliptic
fiber by the order-three automorphism induced by $s\mapsto\zeta_3s$. The three multiplicity-one components
\(D_1,D_2,D_3\) arise from the three fixed points of this action. Hence each
$D_i^*$, $1\leq i\leq3$, can contain the specialization of at most one
fixed torsion point.

The origin gives the zero section meeting \(D_1^*\). If \(3\mid n\), the
other two fixed points are \(n\)-torsion and specialize to \(D_2^*\) and
\(D_3^*\). If \(3\nmid n\), no other fixed torsion point occurs.

The fixed torsion sections therefore specialize to the smooth loci
of the multiplicity-one components. All remaining torsion points lie in orbits of length \(3\) and
specialize to smooth points of \(C_0\). Hence no singular point of
\(\mathcal E_0\) lies on \(T_n\), and no irreducible component of $T_n$ meets
\(C_1^*,C_2^*,C_3^*\).

If $B\subset T_n$ corresponds to an orbit of length $3$,
then $3=B\cdot\mathcal E_0=3(B\cdot C_0)$, so
$B\cdot C_0=1$.
Distinct orbits have distinct images in \(C_0^*\).

When \(3\mid n\), the three fixed points account for three torsion points,
leaving \(n^2-3\) points in length-three orbits. When \(3\nmid n\), only the
origin is fixed, leaving \(n^2-1\) points. Dividing by \(3\) gives the stated
value of \(C_0^*\cdot T_n\), and adding the fixed sections gives the number of
irreducible components.
\end{proof}

\subsubsection*{Type \(III^*\)}

For a fiber of type $III^*$,
Proposition~\ref{prop:IIIstar} gives a fiber
of type $I_0^*$ after a quadratic base change.
The monodromy on the resulting model is
$M_4^2=-I$, so we can use the results
for $I_0^*$.

We label the components of $\mathcal E_0$
so that $C_0,C_1$ have multiplicity one,
with $C_0$ meeting the zero section.
Let $C_3$ be the unique multiplicity-four
component, let $C_2$ be the terminal multiplicity-two
component meeting $C_3$, and $C_4,\ldots,C_7$ the
remaining components. As before, let
$C_i^*$ denote the complement of the
intersection points with the other
components of the central fiber.

\begin{theorem}[Torsion specialization for a type $III^*$ fiber]
\label{thm:type-IIIstar}
With the notation above, the following hold:
\begin{enumerate}
    \item
    \[
    \operatorname{length}(T_n\cap C_0^*)=1,
    \qquad
    \operatorname{length}(T_n\cap C_1^*)=
    \begin{cases}
    1,&2\mid n,\\
    0,&2\nmid n.
    \end{cases}
    \]

    \item No singular point of \(\mathcal E_0\) lies on \(T_n\).

    \item
    \[
    \operatorname{length}(T_n\cap C_2^*)=
    \begin{cases}
    1,&2\mid n,\\
    0,&2\nmid n.
    \end{cases}
    \]
    When it exists, this contribution comes from the unique
    degree-\(2\) component, and
    \[
    B\cdot C_2=1.
    \]

    \item
    \[
    \operatorname{length}(T_n\cap C_3^*)=
    \begin{cases}
    \dfrac{n^2-4}{4},&2\mid n,\\[0.8em]
    \dfrac{n^2-1}{4},&2\nmid n.
    \end{cases}
    \]
    Every component meeting \(C_3^*\) has degree \(4\) over \(\Delta\)
    and satisfies
    \[
    B\cdot C_3=1.
    \]

    \item
    \[
    \operatorname{length}(T_n\cap C_4^*)=\operatorname{length}(T_n\cap C_5^*)=\operatorname{length}(T_n\cap C_6^*)=\operatorname{length}(T_n\cap C_7^*)=
    0.
    \]

    \item
    \[
    \#\operatorname{Irr}(T_n)=
    \begin{cases}
    \dfrac{n^2-4}{4}+3,&2\mid n,\\[0.8em]
    \dfrac{n^2-1}{4}+1,&2\nmid n.
    \end{cases}
    \]
\end{enumerate}
\end{theorem}

\begin{proof}
By Proposition~\ref{prop:finite-monodromy-orbits}, the possible orbit
lengths are \(1,2,\) and \(4\). We apply Theorem~\ref{thm:type-I0star} to the intermediate
$I_0^*$-fiber.

The fixed points of \(M_4\) remain fixed after the quadratic base change.
The origin gives the zero section meeting \(C_0^*\), and, when \(n\) is
even, the second fixed point gives the section meeting \(C_1^*\).

When \(n\) is even, the unique orbit of length \(2\) for \(M_4\) consists
of two points fixed by \(M_4^2=-I\). Thus after the quadratic base change
it becomes two torsion sections on the \(I_0^*\)-model. By
Theorem~\ref{thm:type-I0star}, these meet smooth points of the outer
components. On the original model, they give a single degree-\(2\) component
specializing to \(C_2^*\). Since \(C_2\) has multiplicity \(2\),
\[
2=B\cdot\mathcal E_0=2(B\cdot C_2),
\]
and hence \(B\cdot C_2=1\).

An orbit of length \(4\) for \(M_4\) splits after the quadratic base change
into two orbits of length \(2\) for \(-I\). By
Theorem~\ref{thm:type-I0star}, both specialize to smooth points of the
multiplicity-two component of the intermediate \(I_0^*\)-fiber. On the original model, they give a degree-\(4\) component specializing to the
multiplicity-four component \(C_3\). Thus
\[
4=B\cdot\mathcal E_0=4(B\cdot C_3),
\]
so \(B\cdot C_3=1\).

Theorem~\ref{thm:type-I0star} also shows that all torsion specializations on
the intermediate \(I_0^*\)-fiber occur in smooth loci. Their specialization points on the original model lie in the smooth
loci described above. Hence no singular point of
the \(III^*\)-fiber lies on \(T_n\), and no irreducible component of $T_n$ meets
\(C_4^*,\ldots,C_7^*\).

Finally, when \(n\) is even, two fixed points and the two points in the
length-\(2\) orbit account for four torsion points, leaving
\[
\frac{n^2-4}{4}
\]
orbits of length \(4\). When \(n\) is odd, only the origin is fixed, leaving
\[
\frac{n^2-1}{4}
\]
orbits of length \(4\). Adding the fixed and length-\(2\) components gives
the stated number of irreducible components.
\end{proof}

\subsubsection*{Type \(II^*\)}

For a fiber of type $II^*$,
Proposition~\ref{prop:IIstar} gives a fiber
of type $IV^*$ after a quadratic base change.
The specialization of $n$-torsion on the
resulting model is described by
Theorem~\ref{thm:type-IVstar}.

We label the components of $\mathcal E_0$ so that $C_0$
is the multiplicity-one component met by the
zero section. Let $C_3$ be the unique
multiplicity-six component, $C_1$ the terminal
multiplicity-two component on the arm of
length two, and $C_2$ the terminal
multiplicity-three component meeting $C_3$.
We denote the remaining components by
$C_4,\ldots,C_8$.
As before, let $C_i^*$ denote the complement
of the intersection points with the other
components of the central fiber.

\begin{theorem}[Torsion specialization for a type $II^*$ fiber]
\label{thm:type-IIstar}
With the notation above, the following hold:
\begin{enumerate}
    \item
    \[
    \operatorname{length}(T_n\cap C_0^*)=1.
    \]

    \item No singular point of \(\mathcal E_0\) lies on \(T_n\).

    \item
    \[
    \operatorname{length}(T_n\cap C_1^*)=
    \begin{cases}
    1,&3\mid n,\\
    0,&3\nmid n.
    \end{cases}
    \]
    When it exists, this contribution comes from the unique
    degree-\(2\) component, and
    \[
    B\cdot C_1=1.
    \]

    \item
    \[
   \operatorname{length}(T_n\cap C_2^*)=
    \begin{cases}
    1,&2\mid n,\\
    0,&2\nmid n.
    \end{cases}
    \]
    When it exists, this contribution comes from the unique
    degree-\(3\) component, and
    \[
    B\cdot C_2=1.
    \]

    \item
    \[
    \operatorname{length}(T_n\cap C_3^*)=
    \begin{cases}
    \dfrac{n^2-6}{6},&6\mid n,\\[0.8em]
    \dfrac{n^2-3}{6},
        &3\mid n\text{ and }2\nmid n,\\[0.8em]
    \dfrac{n^2-4}{6},
        &2\mid n\text{ and }3\nmid n,\\[0.8em]
    \dfrac{n^2-1}{6},
        &2\nmid n\text{ and }3\nmid n.
    \end{cases}
    \]
    Every component meeting \(C_3^*\) has degree \(6\) over \(\Delta\)
    and satisfies
    \[
    B\cdot C_3=1.
    \]

    \item
    \[
    \operatorname{length}(T_n\cap C_4^*)=\operatorname{length}(T_n\cap C_5^*)=\operatorname{length}(T_n\cap C_6^*)=\operatorname{length}(T_n\cap C_7^*)=\operatorname{length}(T_n\cap C_8^*)=
    0.
    \]

    \item
    \[
    \#\operatorname{Irr}(T_n)=
    \begin{cases}
    \dfrac{n^2-1}{6}+1,
        &2\nmid n\text{ and }3\nmid n,\\[0.8em]
    \dfrac{n^2-3}{6}+2,
        &3\mid n\text{ and }2\nmid n,\\[0.8em]
    \dfrac{n^2-4}{6}+2,
        &2\mid n\text{ and }3\nmid n,\\[0.8em]
    \dfrac{n^2-6}{6}+3,
        &6\mid n.
    \end{cases}
    \]
\end{enumerate}
\end{theorem}

\begin{proof}
By Proposition~\ref{prop:finite-monodromy-orbits}, the possible orbit
lengths are \(1,2,3,\) and \(6\). We use Theorem~\ref{thm:type-IVstar} on the intermediate
$IV^*$-model.

The origin remains fixed and gives the zero section meeting \(C_0^*\).

If \(3\mid n\), the unique orbit of length \(2\) consists of two points
fixed by \(M_6^2\). After the quadratic base change these become two fixed
torsion sections on the \(IV^*\)-model. By
Theorem~\ref{thm:type-IVstar}, they specialize to two distinct
multiplicity-one components. On the original model, they give a single degree-\(2\)
component specializing to \(C_1^*\). Since \(C_1\) has multiplicity \(2\),
\[
2=B\cdot\mathcal E_0=2(B\cdot C_1),
\]
and hence \(B\cdot C_1=1\).

If \(2\mid n\), the unique orbit of length \(3\) remains an orbit of length
\(3\) after the quadratic base change. By
Theorem~\ref{thm:type-IVstar}, it specializes to the smooth locus of the
multiplicity-three component of the intermediate \(IV^*\)-fiber. On the original model, it gives a degree-\(3\) component specializing to \(C_2^*\). Thus
\[
3=B\cdot\mathcal E_0=3(B\cdot C_2),
\]
so \(B\cdot C_2=1\).

Finally, an orbit of length \(6\) splits after the quadratic base change
into two orbits of length \(3\). Both specialize, by
Theorem~\ref{thm:type-IVstar}, to smooth points of the
multiplicity-three component of the \(IV^*\)-fiber. On the original model, they give
one degree-\(6\) component specializing to the multiplicity-six component
\(C_3\). Hence
\[
6=B\cdot\mathcal E_0=6(B\cdot C_3),
\]
and therefore \(B\cdot C_3=1\).

Theorem~\ref{thm:type-IVstar} shows that all torsion specializations on the
intermediate \(IV^*\)-fiber occur in smooth loci. Their specialization points on the original model lie in the smooth
loci described above. Hence no singular point of
the \(II^*\)-fiber lies on \(T_n\), and no irreducible component of $T_n$ meets
\(C_4^*,\ldots,C_8^*\).

The numbers of degree-\(6\) components are obtained by removing the points
belonging to the shorter orbits and dividing the remaining number by \(6\).
Adding the fixed, degree-\(2\), and degree-\(3\) components when present
gives the stated number of irreducible components.
\end{proof}

\subsection{Irreducibility versus Connectedness}
\label{subsec:additive-connectedness}

The specialization results also determine the connected components of
$T_n$, because distinct irreducible components can meet only on the
central fiber.

For fibers of type $II$, $III$, and $IV$, every non-fixed monodromy
orbit specializes to the unique singular point \(R\). Thus, when $n>1$, their union is connected. The fixed sections
specialize to distinct smooth points and therefore remain separate.

For \(I_0^*,IV^*,III^*\), and \(II^*\), no singular point of the central
fiber lies on \(T_n\), and distinct monodromy orbits specialize to distinct
smooth points. Thus distinct irreducible components of \(T_n\) remain
disjoint.

\begin{corollary}[Connected components for potentially good fibers]
\label{cor:connected-components-additive}
Let \(s(\mathcal E_0,n)\) denote the number of fixed monodromy points on
\(E_t[n]\). For the unstarred fibers,
\[
s(II,n)=1,
\]
\[
s(III,n)=
\begin{cases}
1,&2\nmid n,\\
2,&2\mid n,
\end{cases}
\qquad
s(IV,n)=
\begin{cases}
1,&3\nmid n,\\
3,&3\mid n.
\end{cases}
\]
Then, for \(\mathcal E_0\) of type \(II\), \(III\), or \(IV\),
\[
\#\pi_0(T_n)
=
s(\mathcal E_0,n)+\mathds 1_{n>1}.
\]

For \(\mathcal E_0\) of type
\[
I_0^*,\ IV^*,\ III^*,\ \text{or }II^*,
\]
one has
\[
\#\pi_0(T_n)=\#\operatorname{Irr}(T_n).
\]
\end{corollary}

\begin{proof}
For a fiber of type $II$, $III$, or $IV$, all irreducible components
associated with non-fixed orbits contain the same singular
point \(R\), and hence their union is connected. Such a component exists
exactly when \(n>1\). Each fixed orbit gives a section specializing to a
distinct smooth point, so these sections form separate connected components.

For \(I_0^*,IV^*,III^*\), and \(II^*\), the corresponding specialization
theorems show that distinct irreducible components meet distinct smooth points
of the central fiber and no component passes through a singular point.
Therefore no two irreducible components meet, and the irreducible and
connected components coincide.
\end{proof}

\subsection*{Acknowledgments}

The author would like to express her gratitude
to Fran\c{c}ois Greer for his guidance and support
throughout this project. His many helpful discussions,
valuable suggestions, and feedback have
been essential to the development of this work.

The author used OpenAI's ChatGPT to assist with
English-language editing, LaTeX formatting, and
preliminary consistency checks of the manuscript.
The author takes full responsibility for the
mathematical content.

\end{document}